\documentclass[12pt]{amsart}
 \usepackage{amsfonts,amssymb,amsmath}
 \usepackage{epsfig} 
 \usepackage{graphics} 
\def\l{\langle}
\def\r{\rangle}

\def\p{\partial}

\def\cal{\mathcal}

\newtheorem{Theorem}{Theorem}

\newtheorem{Lemma}{Lemma}[section]
\newtheorem{Proposition}{Proposition}[section]

\theoremstyle{definition}

\newtheorem{remark}{Remark}[section]

\catcode`@=11 \@addtoreset{equation}{chapter}

\catcode`@=12

\title[ Multi-bump solutions on lattices]{\bf Multi-bump solutions of $-\Delta u=K(x)u^{\frac{n+2}{n-2}}$ on Lattices in  ${\mathbb R}^n$}

\author[Y.Y. Li]{Yanyan Li}
\address{\noindent 
School of Mathematical Sciences,
Beijing Normal University,
  Beijing 100875, China
and
Department of Mathematics, Rutgers University, NJ08854,  USA.}
\email{yyli@math.rutgers.edu }

\author[J. Wei]{Juncheng Wei}
\address{\noindent  Department of Mathematics, University of British Columbia, Vancouver, B.C., Canada, V6T 1Z2.
} \email{jcwei@math.ubc.ca}

\author[H. Xu]{Haoyuan Xu}
\address{\noindent  School of Mathematics and Statistics, Huazhong University of Science \& Technology, Wuhan 430074, P. R. China}

\email{hyxu@hust.edu.cn}

\begin{document}

\begin{abstract}
We consider the following semi-linear elliptic equation with
critical exponent:
\[
-\Delta u=K(x)u^{\frac{n+2}{n-2}}, \quad u>0\quad\mbox{in}\quad
{\mathbb R}^n,
\]
where $n\ge 3$, $K>0$ is periodic in $(x_1,..., x_k)$ with $ 1 \leq k <\frac{n-2}{2}$. Under some natural conditions on $K$ near a critical point, we prove the existence of multi-bump solutions where the centers of bumps  can be placed in some lattices in ${\mathbb R}^k$, including infinite lattices. We also show that for $k\ge \frac{n-2}{2}$, no such solutions exist.
\end{abstract}

\date{}
\maketitle

\section{Introduction}\label{S0}
\setcounter{equation}{0}
We consider the following semi-linear elliptic equation with
critical exponent:
\begin{equation}\label{1.1}
-\Delta u=K(x)u^{\frac{n+2}{n-2}}, \quad u>0 \quad\mbox{in}\quad
{\mathbb R}^n.
\end{equation}

\medskip

Associated with (\ref{1.1}) is the following energy  functional
$$
I(u)=\frac{1}{2}\|u\|^2-\frac{n-2}{2n}\int K(x)(u^+)^{\frac{2n}{n-2}},\quad u\in\mathcal D,
$$
where $u^+=\max(u,0)$ and $\mathcal D$ is the Hilbert space defined as the completion of ${\mathcal C}_c^{\infty}({\mathbb R}^n)$ with respect to the scalar product $\l u,v\r=\int_{{\mathbb R}^n}\nabla
u\cdot\nabla v$ and $\|\cdot\|$ denotes the norm of $\mathcal D$. By the maximum principle, a non-zero critical point of $I(u)$ will give rise to a positive solution to equation (\ref{1.1}).

\medskip

When $K\equiv 1$, all solutions of (\ref{1.1}) have been
classified by Caffarelli-Gidas-Spruck \cite{CGS} and are given by:
$$
\sigma_{P,\lambda}(x)=(n(n-2))^{\frac{n-2}{4}}\left(\frac{\lambda}{1+\lambda^2|x-P|^2}\right)^{\frac{n-2}{2}},
$$for any $\lambda>0$ and $P\in {\mathbb R}^n$. See also Obata \cite{O71} and Gidas-Ni-Nirenberg
\cite{GNN} when $u$ has some natural decay as $|x|\to\infty$.

\medskip

When $K$ is positive and periodic, Li proved that (\ref{1.1}) has infinitely many multi-bump solutions for $n\ge 3$ in \cite{L,L1,L2} by gluing approximate solutions into genuine solutions with masses concentrating near isolated sets of maximum points of $K$. Similar results were obtained by Yan in \cite{Yan} if $K(x)$ has a sequence of strict local maximum points tending to infinity. When $K$ is positive and periodic, Xu constructed in
\cite{X1} multi-bump solutions with mass concentrating near critical points of $K$ including saddle points; see also \cite{X2}. When $K(x)$ is a positive radial function with a strict local maximum at $|x|=r_0>0$ and satisfies
$$
K(r)=K(r_0)-c_0|r-r_0|^m+O(|r-r_0|^{m+\theta}),
$$for some constant $c_0>0$, $\theta>0$ and $m\in [2, n-2)$ near $|x|=r_0$, Wei and Yan constructed in \cite{WY} solutions with a large number of bumps concentrating near the sphere $|x|=r_0$ for $n\ge 5$.

\medskip

In this paper, we construct multi-bump solutions of (\ref{1.1}) near critical points of $K(x)$ and the bumps can be placed on {\em arbitrarily many} or even {\em infinitely many} lower dimensional lattice points. Furthermore we show that the dimensional restriction is {\em optimal}.

\medskip

More precisely, we assume the following conditions on $K(x)$:
\begin{enumerate}
\item[{\bf (H1)}] $0<\inf_{{\mathbb R}^n} K\le \sup_{{\mathbb R}^n}K<\infty$;
\item[{\bf (H2)}] $K\in C^1({\mathbb R}^n)$, $K$ is $1$-periodic in its first $k$ variables;

\item[{\bf (H3)}] $0$ is a critical point of $K$ satisfying: there exists some real number $\beta\in (n-2,n)$ such that near $0$,
\begin{equation}\label{1.4}
K(x)=K(0)+\sum_{i=1}^n a_i|x_i|^{\beta}+R(x),
\end{equation}
where $a_i\neq 0$, $\sum_{i=1}^n a_i<0$, and $R(y)$ is
$C^{[\beta]-1,1}$ (up to $[\beta]-1$ derivatives are Lipschitz
functions, $[\beta]$ denotes the integer part of $\beta$) near $0$
and satisfies $\sum_{s=0}^{[\beta]}|\nabla^s
R(y)||y|^{-\beta+s}=o(1)$ as $y$ tends to $0$. Here and
following, $\nabla^s$ denotes all possible partial derivatives of
order $s$.
\end{enumerate}

\medskip

Condition {\bf (H3)} was used by Li in \cite{L2}
for equation (\ref{1.1}). Without loss of generality,  we may assume $K(0)=1$. For any integer $m\ge 1$ and integer $k\in [1,n]$, we define $k$-dimensional lattice
\begin{equation}
\label{qm}
Q_m:=\{\mbox{ all the  integer points in } [0,m]^k\times \{0\}\subset{\mathbb R}^n\},
\end{equation}
  where $0\in {\mathbb R}^{n-k}$. We call $x=(x_1,...,x_n)\in {\mathbb R}^n$ an integer point if all $x_1,...,x_n$ are integers.

\medskip

The main results of this paper can be summarized as follows.

\medskip

\noindent
{\bf Main Theorem:} {\em  For $n \geq 5, 1\leq k <\frac{n-2}{2}$,
there exists $l_0>1$
such that for all $l\ge l_0$ and
all $m\geq 1$ ($m$ can be $+\infty$). There exists a $C^2$ positive solution $u_{Q_{lm}}$ 
with bumps centered close to the lattice set $Q_{lm}$ (defined in (\ref{qm})). 
Furthermore, this is optimal, i.e., if $k \geq \frac{n-2}{2}$, then no such solution exists.}

\medskip

In the following we give more precise statements of the above theorem.

\medskip

For any $\lambda>0$, we define the transformation $S_{\lambda}:u(x)\to u_{\lambda}(x):=\lambda^{-\frac{n-2}{2}}u(\frac{x}{\lambda})$. Then for a solution $u$ of (\ref{1.1}), $u_{\lambda}(x)$ satisfies
\begin{equation}\label{1.2}
-\Delta u_{\lambda}(x)=K_{\lambda}(x)u_{\lambda}(x)^{\frac{n+2}{n-2}},\quad u_{\lambda}>0,\quad x\in{\mathbb R}^n.
\end{equation}
Since $K$ is 1-periodic in its first $k$ variables, $K_{\lambda}(x):=K(\frac{x}{\lambda})$ is $\lambda$-periodic in its first $k$ variables.

\medskip

For any positive integer $l$, let
\begin{equation}\label{1.2a}
\lambda=l^{\frac{n-2}{\beta-n+2}},
\end{equation}where $n-2<\beta<n$. (Throughout this paper, $l$ and $\lambda$ will satisfy the relation (\ref{1.2a}).)

\medskip

We scale the lattice $Q_m$ as
\begin{equation}
\label{1.2ab}
X_{l,m}=\{\lambda l x|x\in Q_m\}.
\end{equation}

\medskip

For convenience, we order the points in $X_{l,m}$ as $\{X^i\}_{i=1}^{(m+1)^k}$.

\medskip

For $i=1,..., (m+1)^k$, let $P^i\in B_{\frac{1}{2}}(X^i)$ and $\Lambda_i>0$. We use notations $P:=\{P^i\}_{i=1}^{(m+1)^k}$ and $\Lambda:=\{\Lambda_i\}_{i=1}^{(m+1)^k}$. Then
$$
W_m(x,P,\Lambda)=\sum_{i=1}^{(m+1)^k} \sigma_{P^i,\Lambda_i}(x)
$$
is an approximate solution of (\ref{1.2}) when $l$ is much larger than $\max_{i} \Lambda_i $. When there is no confusion, we denote
$W_m(x,P,\Lambda)$ by $W_m(x)$.

\medskip

For a fixed lattice $X_{l,m}$ and $\tau>1$, and for functions $\phi, f\in L^{\infty}({\mathbb R}^n)$, let
\begin{equation}\label{2.1}
\|\phi\|_*=\sup_{y\in{\mathbb R}^n}\left(\gamma(y)\sum_{i=1}^{(m+1)^k}\frac{1}{(1+|y-X^i|)^{\frac{n-2}{2}+\tau}}\right)^{-1}|\phi (y)|,
\end{equation}
and
\begin{equation}\label{2.2}
\|f\ |_{**}=\sup_{y\in{\mathbb R}^n}\left(\gamma(y)\sum_{i=1}^{(m+1)^k}\frac{1}{(1+|y-X^i|)^{\frac{n+2}{2}+\tau}}\right)^{-1}|f(y)|,
\end{equation}
where
\begin{equation}
\gamma(y)=\min\left(1, \min_{i=1}^{(m+1)^k} (\frac{1+|y-X^i|}{\lambda})^{\tau-1}\right).
\end{equation}

The above weighted norms depend on $\tau$ and the lattice $X_{l,m}$, since $\tau$ and $X_{l,m}$ are chosen. When there is no confusion, we just denote them as above.

\medskip

\begin{remark}
\em
 Without the function $\gamma(y)$, the above weighted norms are used by Wei-Yan in \cite{WY}. Similar weighted norms can be found in \cite{DFM,RW1,RW2}. The reason for introducing $\gamma(y)$ in the definition of the weighted norms is crucial in our proofs. We shall comment more on this later.
\end{remark}

We now state the first result.

\begin{Theorem}\label{T2}
For  $n\ge 5$, $1\le k<\frac{n-2}{2}$ and $m\ge 1$, assume that $K$ satisfies conditions {\bf (H1)}, {\bf (H2)} and {\bf (H3)}. Then there exists a $\tau_0(n,k)\in (k,\frac{n-2}{2}]$, such that for any $\tau$ satisfying $k\le\tau<\tau_0$, there exist positive constants $C_1$, $C_2$, $C$ and an integer $l_0$ depending only on $K$, $n$, $\beta$, $\tau$, such that for any integer $l\ge l_0$ and $\lambda=l^{\frac{n-2}{\beta-n+2}}$, equation (\ref{1.1}) has a positive $C^2$ solution $u$ satisfying
$$
\|S_{\lambda} u-W_m\|_*\le  C\lambda^{\tau-\frac{n+2}{2}},
$$with $C_1< \Lambda_i< C_2$ for all $i$ and $\max_{1\le i\le (m+1)^k}|P^i-X^i|\to 0$ as $l\to\infty$, uniformly in $m$.
\end{Theorem}

If we allow the estimates to depend on $m$ (e.g., the size of $l_0$), then $m$-bump solutions have been constructed by Xu in \cite{X1} for every $1\le k\le n$ under the same assumption on $K$, see also \cite{X2}. As mentioned earlier, Li constructed in \cite{L, L1, L2} such $m$-bump solutions near isolated sets of maximum points. The ansatz used in \cite{L,L1,L2} is a variational method as in \cite{CP,CP1} of Coti Zelati-Rabinowitz and \cite{SE} of S\'{e}r\'{e} which glues approximate solutions into genuine solutions. On the other hand, the method used in \cite{X1,X2} is gluing via implicit function theorem (or a nonlinear Lyapunov-Schmidt technique), the same as that in our proof of Theorem \ref{T2}. Such Lyapunov-Schmidt reduction methods have been developed and used by many aithors. We shall make some comments at the end of this section.

\medskip

The {\em novelty} and the {\em main difficulty} in the proof of Theorem \ref{T2} is that all the estimates are {\em independent of $m$} and the results {\em are optimal} (see Theorem 3 below). Thus we may construct multi-bump solutions on {\em an infinite lattice} by letting $m\to\infty$ while keeping $l$ fixed (see Theorem 2 below). The new $m$-independent estimates are obtained by using the new weighted norm $\|\cdot\|_*$  (defined at (\ref{2.1})) as compared to $\|\cdot\|$ of $\mathcal D$ used in \cite{X1}. Roughly speaking, $\|\cdot\|$ norm adds up errors near each bump, while $\|\cdot\|_*$ norm measures maximum of errors near each bump.
The reason for introducing $\gamma(y)$ is to {\em localize} the estimate and to obtain better decay estimates near each bump. This is crucially needed when we deal with the higher dimensional lattice case in which $ 2\leq k\leq \tau$.
As we mentioned earlier,  Wei-Yan \cite{WY} used a similar norm in which $\gamma (y)\equiv 1$. They required  that the number $ \tau$ must be $1+\bar{\eta}$ with $\bar{\eta}>0$ being small. Thus the norms used in \cite{WY} are only suitable for concentration on one dimensional set (like circles as in \cite{WY}). Our norms work for any higher dimensional concentration as long as the dimension $ k <\frac{n-2}{2}$ (which is optimal). This is one of major technical advances in this paper.

\medskip

Let
$$
{\mathbb Z}^k:=\{\mbox{all integer points in } {\mathbb R^k}\times \{0\}\},\mbox{ where } 0\in {\mathbb R}^{n-k},
$$and, for an integer $i\in [0,k]$,
$$
{\mathbb R}^k_i\times\{0\}:=\{(x_1,...,x_k,0,...,0)\in {\mathbb R}^n|x_1,...,x_i\ge 0\}.
$$

\medskip

We often write ${\mathbb R}^k_i\times\{0\}$ as ${\mathbb R}^k_i$ when there is no confusion. Note also that ${\mathbb R}^k_0={\mathbb R}^k$.
Consider infinite lattices
$$
Y^i\equiv Y^{k,i}:={\mathbb Z}^k\cap {\mathbb R}^k_i
$$ and their scaled versions
$$
X_l^i=\lambda lY^i.
$$ Define
$$
W_l^i=\sum_{X\in X_l^i}\sigma_{P(X),\Lambda(X)},
$$ where $P(X)\in B_{\frac{1}{2}}(X)$ and $\Lambda(X)\in (C_1,C_2)$.

\medskip

\begin{Theorem}\label{T1}
For $n\ge 5$, $1\le k<\frac{n-2}{2}$ and $0\le i\le k$, assume that $K$ satisfies conditions {\bf (H1)}, {\bf (H2)} and {\bf (H3)}. Let $\tau$, $C_1$, $C_2$, $C$ and  $l_0$ be as in Theorem \ref{T2}. Then for any integer $l\ge l_0$ and $\lambda=l^{\frac{n-2}{\beta-n+2}}$, equation (\ref{1.1}) has a positive $C^2$ solution $u$ satisfying
$$
\|S_{\lambda} u-W_l^i\|_*\le  C\lambda^{\tau-\frac{n+2}{2}},
$$ with $\Lambda(X)\in (C_1, C_2)$ for all $X\in X_l^i$ and $|P(X)-X|\to 0$ as $l\to\infty$ uniformly in $X$.
\end{Theorem}

Solutions $u$ constructed in Theorem \ref{T1} have infinitely many bumps. Indeed, $u$ is close to $\sigma_{P(X),\Lambda(X)}$ near every lattice point $X\in X_l^i$.

\medskip

Theorem \ref{T1} follows from Theorem \ref{T2} by a limiting argument as follows. For $l\ge l_0$, let $X_l^i$ be an infinite lattice as in Theorem \ref{T1}. By Theorem \ref{T2}, we have solutions $u_m$ of (\ref{1.1}) for all $m$. For each $m$, we can find $x_m\in X_{l,m}$, such that $X_{l,m}-x_m$ is monotonically increasing in $m$ and
$$
\cup_{m=1}^{\infty}(X_{l,m}-x_m)=X_l^i.
$$

Let
$$
(S_{\lambda}\hat{u}_m)(x)=(S_{\lambda}u_m)(x+x_m).
$$

Then $\hat{u}_m$ satisfies the same equation as $u_m$ due to the periodicity of $K$. See section \ref{S2} for details.

\medskip

\begin{remark}
\em
Theorem \ref{T1} shows a new phenomena that infinite-bump solutions can be constructed for semilinear elliptic equations with critical exponents. For subcritical exponent semilinear elliptic equations, infinite-bump solutions were constructed by Coti Zelati and Rabinowitz in \cite{CP,CP1} and S\'{e}r\'{e} in \cite{SE}. There are two main differences (and difficulties) between the subcritical exponent problem (treated in \cite{CP, CP1}) and critical exponent problem: first, there is an extra loss of compactness-the scaling invariance parameter. Second, there is the difficulty of controlling  the algebraic decaying in an infinite lattice setting.  (In \cite{CP, CP1}, the decay rate is exponential.) As far as we know, this paper seems to be the first in obtaining the existence of solutions for critical exponent problems with {\em infinitely many bumps}.

\end{remark}

\medskip

In an unpublished note \cite{L7}, Li showed that the conclusion of Theorem \ref{T2} and Theorem \ref{T1} are false if $n\ge 3$ and $k\ge n-2$. More specifically, let $K$ be a positive $C^1$ function which is periodic in each variable, satisfying, for some constant $\beta>n-2$, $(*)_{\beta}$ condition for some positive constants $L_1$ and $L_2$ in ${\mathbb R}^n$:
$$|\nabla K_i|\le L_1,\quad \mbox{in}\quad {\mathbb R}^n,$$ and, if $\beta\ge 2$, that $K\in C^{[\beta]-1,1}_{loc}({\mathbb R}^n)$,

$$
|\nabla^s K_i(y)|\le L_2|\nabla K_i(y)|^{\frac{\beta-s}{\beta-1}},\quad\mbox{for all}\quad  2\le s\le [\beta],\quad \forall y\in {\mathbb R}^n.
$$  Then for $n\ge 3$ and $k\ge n-2$, there is no $C^2$ solution of (\ref{1.1}) satisfying, for some $R$, $\epsilon>0$ and $0\le i\le k$,
\begin{equation}\label{A2}
\inf_{x\in {\mathbb R}_i^k}\sup_{B_R(x)}u\ge \epsilon.
\end{equation}

\medskip

$(*)_{\beta}$ condition was introduced in \cite{L2}. If a positive $C^1$ periodic function $K$ is of the form (\ref{1.4}) near every critical point of $K$, then $K$ satisfies $(*)_{\beta}$ in ${\mathbb R}^n$. Also $(*)_{\beta_1}$ implies $(*)_{\beta_2}$ if $\beta_1\ge \beta_2$. If a function $K$ in (\ref{1.1}) satisfies $(*)_{\beta}$ for some $\beta>n-2$, then solutions in any bounded region can only have isolated simple blow up points, see \cite{L2}.

\medskip

Our next theorem improves this result to cover $k\ge \frac{n-2}{2}$, which is optimal in view of Theorem \ref{T2} and Theorem \ref{T1}.
\begin{Theorem}\label{T4}
For $n\ge 3$ and $k\ge\frac{n-2}{2}$, let $K$ be as above. Then there is no $C^2$ solution of (\ref{1.1}) satisfying (\ref{A2}) for some $R$, $\epsilon>0$ and $0\le i\le k$.
\end{Theorem}

\medskip

\begin{remark}
\em
The solutions constructed in Theorem \ref{T1} satisfy (\ref{A2}) for some $R$, $\epsilon>0$ and $0\le i\le k$. Since the hypotheses on $K$ in Theorem \ref{T4} are stronger than that in Theorem \ref{T1}, the assumption $k<\frac{n-2}{2}$ in Theorem \ref{T1} is optimal.
\end{remark}

\medskip

We  end the introduction with   some remarks and history on the finite/infinite dimensional reduction method. The original finite dimensional Liapunov-Schmidt reduction method was first introduced in a seminal paper by Floer and Weinstein \cite{FW} in their construction of single bump solutions to one dimensional  nonlinear Schrodinger equations (Oh \cite{Oh} generalized to high dimensional case). On the other hand, Bahri \cite{B} and Bahri-Coron \cite{BC1} developed the reduction method for critical exponent problems. In the last fifteen years,  there are  renewed efforts in refining the finite dimensional reduction method by many authors. When combined with variational methods, this reduction becomes "localized energy method". For subcritical exponent problems, we refer to Ambrosetti-Malchiodi \cite{AM},  Gui-Wei \cite{GW}, Malchiodi \cite{ma}, Li-Nirenberg \cite{L3},  Lin-Ni-Wei \cite{LNW}, Ao-Wei-Zeng \cite{AWZ} and the references therein. The localized energy method in degenerate setting is done by Byeon-Tanaka \cite{BT1, BT2}. For critical exponents,  we refer to   Bahri-Li-Rey \cite{BLR}, Del Pino-Felmer-Musso \cite{DFM}, Rey-Wei \cite{RW1, RW2} and Wei-Yan \cite{WY} and the references therein.   Many new features of the finite dimensional reduction are found. Our current work contributes to this part of reduction method and gives {\em an optimal treatment} for critical exponent problems. In recent years, there are new interests in extending the finite dimensional reduction method to treat high dimensional concentration phenomena.  This is the infinite dimensional reduction method  and has become  very useful in constructing high dimensional concentration solutions. For compact manifold case, we refer to Del Pino-Kowalczyk-Wei \cite{DKW-CPAM, DKW-JDG} and Pacard-Ritore \cite{Pacard}, and for noncompact manifolds, we refer to Del Pino-Kowalczyk-Pacard-Wei  \cite{DKPW}, Del Pino-Kowalczyk-Wei \cite{DKW-DG}  and the references therein. A notable application of this infinite dimensional reduction method is the construction of  counterexample to De Giorgi' s Conjecture in large dimensions by M. del Pino, M. Kowalczyk and Wei \cite{DKW-DG}.

\medskip

Throughout the paper, we will use the superscript to stand for a sequence of points in ${\mathbb R}^n$ and the subscript to stand for numbers or the coordinates of a point in ${\mathbb R}^n$ unless otherwise stated.

\medskip

The paper is organized as follows. In section \ref{S1}, we carry out the Lyapunov-Schmidt reduction. In section \ref{S2}, we solve the finite dimensional problems and prove Theorem \ref{T2} and Theorem \ref{T1}. In section \ref{S3}, we prove Theorem \ref{T4}. In Appendices, we prove some basic lemmas that will be used throughout the paper.

\bigskip

\noindent
{\bf Acknowledgments:} The research of the first author is partially 
supported by an NSF grant, the research of the second author is
partially supported by a General Research Fund from RGC of Hong Kong.

\section{Finite Dimensional Reduction}\label{S1}
\setcounter{equation}{0}

In this section, we perform a finite-dimensional reduction to the re-scaled equation (\ref{1.2}).  Let $\lambda $ and $l$ be given at (\ref{1.2a}) and $ K_\lambda (x)= K(\frac{x}{\lambda})$.

\medskip

Consider the following energy functional
$$
I_{\lambda}(u)=\frac{1}{2}\|u\|^2-\frac{n-2}{2n}\int K_{\lambda}(u^+)^{\frac{2n}{n-2}},\quad u\in\mathcal D.
$$
Then any nonzero critical point of $I_{\lambda}$ gives rise to a solution to (\ref{1.2}). For $X^i\in X_{l,m}$, let $P^i\in B_{\frac{1}{2}}(X^i)$. For any positive constants $C_1<C_2$, let $\Lambda_i\in (C_1,C_2)$ and denote $\sigma_i=\sigma_{P^i,\Lambda_i}$. For a small $\rho>0$, let
$$
\Sigma=\{\sum_i (1+\epsilon_i)\sigma_{P_i,\Lambda_i}| P^i\in B_{\frac{1}{2}}(X^i),\Lambda_i\in (C_1,C_2),|\epsilon_i|<\rho\}.
$$

\medskip

Then $\Sigma$ is a smooth $(n+2)(m+1)^k$ dimensional sub-manifold in $\mathcal D$. When $l$ is large enough, according to Proposition 2 of Bahri-Coron \cite{BC}, every function in a small tubular neighborhood of $\Sigma$ in $\mathcal D$ can be uniquely parameterized as
$$
u(x)=\sum_i(1+\epsilon_i)\sigma_i(x)+\phi(x):=\bar{W}_m(x)+\phi(x),
$$ where $\phi(x)$ is the unique minimizer of
$$
\min_{w\in \Sigma }\|u-w\|.
$$

\medskip

In particular, $\phi\in{\mathcal E}$, a subspace of $\mathcal D$ defined as

$$
\mathcal E:=\left\{\phi\in\mathcal D |\langle Z_{i,j},\phi\rangle=0, \langle\sigma_i,\phi\rangle=0, i=1,...,(m+1)^k, j=1,...,n+1 \right\},
$$ with $Z_{i,j}$ defined as

$$
Z_{i,j}=\frac{\partial \sigma_{P^i,\Lambda_i}}{\partial {P^i}_j}\quad\mbox{for $1\le j\le n$ and}\quad Z_{i,n+1}=\frac{\partial\sigma_{P^i,\Lambda_i}}{\partial \Lambda_i},
$$where ${P^i}_j$ means the $j$-th coordinate of $P^i$. The size of the tubular neighborhood is independent of $l$ and $m$. Obviously, $\mathcal E$ depends on the choice of $P$ and $\Lambda$.

\medskip

In the form $u=\bar{W}_m+\phi$, $I_{\lambda}(u)=I_{\lambda}(\bar{W}_m+\phi)$ can be written as
\begin{equation}\label{2.9}
J(P, \epsilon, \Lambda,\phi)=\frac{1}{2}\int_{{\mathbb R}^n}|\nabla (\bar{W}_m+\phi)|^2-\frac{n-2}{2n}\int_{{\mathbb R}^n}K_{\lambda}(y)\left((\bar{W}_m+\phi)^+\right)^{\frac{2n}{n-2}},
\end{equation}where $\bar{W}_m\in\Sigma$ and $\phi\in\mathcal E$.

\medskip

By definition, $u=\bar{W}_m+\phi$ is a critical point of $I_{\lambda}$ if and only if
\begin{equation}\label{2.9a}
\frac{\partial J}{\partial \phi}(P,\epsilon,\Lambda,\phi)=0,
\end{equation}

\begin{equation}\label{2.9b}
\frac{\partial J}{\partial \epsilon}(P,\epsilon,\Lambda,\phi)=0,
\end{equation}

\begin{equation}\label{2.9c}
\frac{\partial J}{\partial P}(P,\epsilon,\Lambda,\phi)=0,
\end{equation}
\begin{equation}\label{2.9d}
\frac{\partial J}{\partial \Lambda}(P,\epsilon,\Lambda,\phi)=0.
\end{equation}

\medskip

We will use Lyapunov-Schmidt reduction method to solve these equations. More specifically, for fixed $P$ and $\Lambda$, we first solve (\ref{2.9a}) and (\ref{2.9b}), finding solutions $\phi(P,\Lambda)$ and $\epsilon(P,\Lambda)$ which are $C^1$ in $P$ and $\Lambda$. Then we use the Brouwer fixed point theorem to solve the finite dimensional problems (\ref{2.9c}) and (\ref{2.9d}).

\medskip

Let
\begin{equation}\label{2.10.1}
F(\phi,\epsilon):=\frac{\partial J}{\partial \phi}(P,\epsilon,\Lambda,\phi);
\end{equation}
\begin{equation}\label{2.10.2}
G_i(\phi,\epsilon):=\frac{\partial J}{\partial \epsilon_i}(P,\epsilon,\Lambda,\phi), \ i=1,..., (m+1)^k.
\end{equation}

\medskip

The explicit expressions for $F$ and $G$ are as follows
$$
G_i(\phi,\epsilon)=\sum_j(1+\epsilon_j)\l\sigma_i,\sigma_j\r-\int_{{\mathbb R}^n}K_{\lambda}\left((\bar{W}_m+\phi)^+\right)^{\frac{n+2}{n-2}}\sigma_i,
$$
$$
F(\phi,\epsilon)=\phi-P_{\mathcal{E}}(-\Delta)^{-1}K_{\lambda}\left((\bar{W}_m+\phi)^+\right)^{\frac{n+2}{n-2}},
$$
where $P_{\mathcal{E}}$ is the orthogonal projection of $\mathcal D$ onto $\mathcal E$.

\medskip

For fixed $(P,\Lambda)$, setting
\begin{equation}\label{2.10.3}
N(\phi,\epsilon)=(F,G)(\phi,\epsilon):\mathcal{E}\times {\mathbb R}^{(m+1)^k}\to\mathcal{E}\times {\mathbb R}^{(m+1)^k}
\end{equation}
 where $G=(G_1,...,G_{(m+1)^k})$. With the aid of the implicit function theorem, we solve
\begin{equation}
N(\phi,\epsilon)=0.
\end{equation}
We shall, as in \cite{L3}, make use of the following form of the implicit function theorem:
\begin{Lemma}\label{L1.1}
\ (Brezis and Nirenberg)\
Let $\mathcal{X}$, $\mathcal Y$ be Banach spaces, $a>0$, $B_a=B_a(z_0)=\{z\in
{\cal X}:\|z-z_0\|\le a\}$. Suppose that $F$ is a $C^1$ map of
$B_a$ into $\cal Y$,with $F'(z_0)$ invertible, and satisfying, for
some $0<\theta<1$,
$$
\|F'(z_0)^{-1}F(z_0)\|\le (1-\theta)a,
$$

$$
\|F'(z_0)^{-1}\|\|F'(z)-F'(z_0)\|\le \theta\quad\forall z\in B_a.
$$
Then there is a unique solution in $B_a$ of $F(z)=0$.
\end{Lemma}

Define
$$
\mathcal{M}=\{u\in L^{\infty}({\mathbb R}^n)|\quad\|u\|_*<\infty\},
$$

$$
\mathcal{\tilde D}=\{u\in L^{\infty}({\mathbb R}^n)|\quad \|u\|_{**}<\infty\}.
$$

\medskip

We will solve the equation $N(\phi,\epsilon)=0$ under the weak sense in Banach spaces with norms related to $\|\cdot\|_*$, defined at (\ref{2.1}). Since $\sigma_i$, $Z_{i,j}\in \mathcal{M}$, we define a subspace of $\mathcal M$ by
\begin{equation}\label{l1.1.1}
\begin{array}{l}
\mathcal{\tilde M}:=\{\phi\in \mathcal{M}|\int_{{\mathbb R}^n}\phi\sigma_i^{\frac{n+2}{n-2}}=0,\quad \int_{{\mathbb R}^n}\phi\sigma_i^{\frac{4}{n-2}}Z_{i,j}=0,\\
\\
\mbox{for all}\quad i=1,...,(m+1)^k,\quad j=1,...,n+1\}.
\end{array}
\end{equation}
For functions $\phi\in\mathcal M$, $\l \phi,\sigma_i\r$ (or $\l \sigma_i,\phi\r$) should be understood as $-\int_{{\mathbb R}^n}\phi\Delta\sigma_i$ (Similar statement also works for $Z_{i,j}$).
Without introducing new symbols, we still use $P_{\mathcal E}$ to denote the orthogonal projection from $\mathcal{M}\to \mathcal{\tilde M}$, which can be defined as follows.

\medskip

Let $\{f_1,...,f_{(n+2)(m+1)^k}\}$ be an orthonormal basis of the $span \{\{\sigma_i\},\{Z_{i,j}\}\}$ obtained by Gram-Schmidt procedure. Then, for $\phi\in \mathcal D$,
$$
\begin{array}{rl}
P_{\mathcal E}\phi=&\phi-\sum_{i=1}^{(n+2)(m+1)^k}\l \phi,f_i\r f_i \\
\\
=&\phi+\sum_{i=1}^{(n+2)(m+1)^k}\left(\int_{{\mathbb R}^n}\phi\Delta f_i\right)f_i.
\end{array}
$$

Therefore, for any $\phi\in \mathcal M$, we define
$$
P_{\mathcal E}\phi=\phi+\sum_i\left(\int_{{\mathbb R}^n}\phi\Delta f_i\right)f_i.
$$

Then $P_{\mathcal E}\phi\in \mathcal{\tilde{M}}$ for every $\phi\in\mathcal M$. It is clear that if $\phi\in{\mathcal D}\cap {\mathcal M}$, then the definitions are the same.

\medskip

For any $h\in \mathcal{\tilde D}$, we define $(-\Delta)^{-1} h$ as
\begin{equation}\label{l1.1.2}
(-\Delta)^{-1} h:=\frac{1}{n(n-2)\omega_n}\int_{{\mathbb R}^n}\frac{h(z)}{|y-z|^{n-2}}dz
\end{equation}where $\omega_n$ is the volume of a unit ball in ${\mathbb R}^n$. It is not difficult to see that (\ref{l1.1.2}) makes sense. Lemma \ref{La.11} in appendix A shows that $(-\Delta)^{-1} h\in \mathcal M$ and the proof of Lemma \ref{L3} actually shows that $P_{\mathcal{E}}(-\Delta)^{-1}$ is a well defined map from $\mathcal {M}\to \mathcal{\tilde M}$. Therefore, we can view $N(\phi,\epsilon)$ as a map from $\mathcal{\tilde{M}}\times {\mathbb R}^{(m+1)^k}\to\mathcal{{\tilde M}}\times {\mathbb R}^{(m+1)^k}$. When $(\phi,\epsilon)$ solves $N(\phi,\epsilon)=0$ in $\mathcal{\tilde{M}}\times {\mathbb R}^{(m+1)^k}$, then $(\phi,\epsilon)$ is automatically a solution to (\ref{2.10.3}) in the weak sense in the original spaces.

\medskip

We will apply Lemma \ref{L1.1} to $N$ at $(0,0)$ in the Banach spaces $\mathcal{X}=\mathcal{Y}=\mathcal{\tilde{M}}\times {\mathbb R}^{(m+1)^k}$, with norm
$$
\|(\phi,\epsilon)\|=\max(\|\phi\|_*,\lambda^{\tau-1}|\epsilon|),
$$where $|\epsilon|=\max_i |\epsilon_i|$.

\medskip

We have the following proposition.

\begin{Proposition}\label{p2}
Under the assumptions of Theorem \ref{T2}, when $l\ge l_0$, equation (\ref{2.10.3}) has a unique solution $\phi(P,\Lambda)$, $\epsilon(P,\Lambda)$ in $\mathcal{{\tilde M}}\times {\mathbb R}^{(m+1)^k}$, with
$$
\|(\phi,\epsilon)\|\le \frac{C}{\lambda^{\frac{n+2}{2}-\tau}}.
$$Furthermore $\phi(P,\Lambda)$ and $\epsilon(P,\Lambda)$ are $C^1$ in $P$ and $\Lambda$.
\end{Proposition}
\begin{proof}
The proof will be carried out in several steps.

\medskip

\noindent
{\bf Step 1.} We first show that, for some constant $C>0$, independent of $m$ and $l$, we have the following crucial estimate for the error
\begin{equation}\label{p1.1}
\|N(0,0)\|\le \frac{C}{\lambda^{\frac{n+2}{2}-\tau}}.
\end{equation}

Observe that
$$
F(0,0)=-P_{\mathcal{E}}(-\Delta)^{-1}K_{\lambda}W_m^{\frac{n+2}{n-2}}.
$$

\medskip

In view of Lemma \ref{La.11} and the fact that $\l F(0,0),\sigma_i\r=0$ for all $i$, we just need to estimate  $\|l_m\|_{**}$, where the error becomes
\begin{equation}
\label{error}
l_m=K_{\lambda}W_m^{\frac{n+2}{n-2}}-\sum_i\sigma_i^{\frac{n+2}{n-2}}.
\end{equation}

\medskip

 Denote $\Omega_i=\{y\in {\mathbb R}^n,  |y-X^i|\le |y-X^j|, \quad\mbox{for all}\quad j\neq i\}$, $B_i=B_{\lambda l}(X^i)$, and $B_{i,m}=B_{\max(\frac{m}{4},1)\lambda l}(X^i)$. Certainly $ {\mathbb R}^n= \cup_{i} \Omega_i$.

\medskip

For $y\in\Omega_i$, it holds that
$$
|K(\frac{y}{\lambda})W_m^{\frac{n+2}{n-2}}-\sum_i\sigma_i^{\frac{n+2}{n-2}}|\le C|K(\frac{y}{\lambda})-1|\sigma_i^{\frac{n+2}{n-2}}+C\left(\hat{W}_{m,i}^{\frac{n+2}{n-2}}+\sigma_i^{\frac{4}{n-2}}\hat{W}_{m,i}\right)+C\sum_{j\neq i}\sigma_j^{\frac{n+2}{n-2}},
$$where $\hat{W}_{m,i}=\sum_{j\neq i}\sigma_j$. We apply  Lemma \ref{L5} to estimate each term on the right hand side of the above inequality. First we estimate $ \hat{W}_{m, i}^{\frac{n+2}{n-2}}$.

\medskip

By Lemma \ref{L5}, if $y\in\Omega_i\cap B_i^c \cap B_{i,m}$, we have
$$\begin{array}{c}
\hat{W}_{m,i}^{\frac{n+2}{n-2}}\le \frac{C}{(\lambda l)^{\frac{4}{n-2}k}}\sum_j\frac{1}{(1+|y-x^j|)^{n+2-\frac{4}{n-2}k}}\\
\\
\le C\sum_j\frac{1}{(1+|y-x^j|)^{\frac{n+2}{2}+\tau}}\frac{1}{(\lambda l)^{\frac{n+2}{2}-\tau}}.
\end{array}
$$

\medskip

For $y\in\Omega_i\cap B_i$ we obtain
$$
\hat{W}_{m,i}^{\frac{n+2}{n-2}}\le \frac{C}{(\lambda l)^{n+2}}\le\left\{\begin{array}{c}
\frac{C}{(\lambda l)^{\frac{n+2}{2}-\tau}}\sum_j\frac{1}{(1+|y-x^j|)^{\frac{n+2}{2}+\tau}},\quad y\in\Omega_i\cap B_{\lambda}(x^i)^c\\
\\
\left(\frac{1+|y-x^i|}{\lambda}\right)^{\tau-1}\sum_j\frac{1}{(1+|y-x^j|)^{\frac{n+2}{2}+\tau}}\frac{\lambda^{\tau-1}}{(\lambda l)^{\frac{n}{2}}},\quad y\in B_{\lambda}(x^i)\cap\Omega_i.
\end{array}\right.
$$

If $y\in B_{i,m}^c\cap\Omega_i$,  by Lemma \ref{L5}, we have first
$$
\hat{W}_{m, i}^{\frac{n+2}{n-2}} \leq  \left(\sum_j\frac{1}{(1+|y-x^j|)^{n-2}}\right)^{\frac{n+2}{n-2}}
$$
\begin{equation}
\label{w36}
\le C\frac{1}{(1+|y-x^i|)^{n+2}}m^{\frac{n+2}{n-2}k}.
\end{equation}

On the other hand
\begin{equation}
\label{w34}
\sum_j\frac{1}{(1+|y-x^j|)^{n}}\ge C\frac{1}{(1+|y-x^i|)^{n}}(1+C[\frac{m}{2}])^k.
\end{equation}

It is not hard to see that when $y\in B_{i,m}^c\cap\Omega_i$, there is a constant $C$ independent of $m$ such that for $\frac{n-2}{2} > k$,
\begin{equation}
\label{w35}
\frac{1}{(1+|y-x^i|)^{2}}m^{\frac{n+2}{n-2}k}\le C\frac{1}{(\lambda l)^{\frac{4}{n-2}k}}(1+C[\frac{m}{2}])^k.
\end{equation}

From (\ref{w36}), (\ref{w34}) and (\ref{w35}), we conclude that for $y \in B_{i, m}^c \cap \Omega_i$
$$
\hat{W}_{m, i}^{\frac{n+2}{n-2}} \le C\sum_j\frac{1}{(1+|y-x^j|)^{\frac{n+2}{2}+\tau}}\frac{1}{(\lambda l)^{\frac{4}{n-2}k}}.
$$

\medskip

Combining the previous estimates together, we get
\begin{equation}
\label{hatWest}
\|\hat{W}_{m,i}^{\frac{n+2}{n-2}}\|_{**}\le \frac{C}{(\lambda l)^{\frac{n+2}{2}-\tau}}.
\end{equation}

 \medskip

Similarly, we have
$$
\|\sigma_i^{\frac{4}{n-2}}\hat{W}_{m,i}\|_{**},\quad \|\sum_{j\neq i}\sigma_i^{\frac{n+2}{n-2}}\|_{**}\le \frac{C}{(\lambda l)^{\frac{n+2}{2}-\tau}}.
$$

\medskip

Now if $y\in\Omega_i$ such that $|y-x^i|\le\lambda$, we get
$$\begin{array}{c}
|K(\frac{y}{\lambda})-1|\sigma_i^{\frac{n+2}{n-2}}\le \frac{C|y-x^i|^{\beta}}{\lambda^{\beta}}\frac{1}{(1+|y-x^i|)^{n+2}}\\
\\
\le \frac{C}{\lambda^{\frac{n+2}{2}-\tau}}\gamma(y)\sum_j\frac{1}{(1+|y-x^j|)^{\frac{n+2}{2}+\tau}}.
\end{array}
$$

\medskip

On the other hand,  if $y\in\Omega_i$ such that $|y-x^i|\ge\lambda$,
$$
|K(\frac{y}{\lambda})-1|\sigma_i^{\frac{n+2}{n-2}}\le \frac{C}{(1+|y-x^i|)^{n+2}}\le \frac{C}{\lambda^{\frac{n+2}{2}-\tau}}\sum_j\frac{1}{(1+|y-x^j|)^{\frac{n+2}{2}+\tau}}.
$$

Thus we obtain the estimate for $\| (K_\lambda -1) \sigma_i^{\frac{n+2}{n-2}} \|_{**}$.

\medskip

Combining the above inequalities, we get
$$
\|l_m\|_{**}\le \frac{C}{\lambda^{\frac{n+2}{2}-\tau}},
$$
and therefore by Lemma \ref{La.11},
$$
\|F(0,0)\|_*\le \frac{C}{\lambda^{\frac{n+2}{2}-\tau}}.
$$

\medskip

 Next we estimate $ G(0,0)$. For each $i$,
$$
G_i(0,0)=\sum_j\l\sigma_i,\sigma_j\r-\int_{{\mathbb R}^n}K_{\lambda}W_m^{\frac{n+2}{n-2}}\sigma_i.
$$

\medskip

It is easy to see that
$$
W_m^{\frac{n+2}{n-2}}=\sigma_i^{\frac{n+2}{n-2}}+C\sigma_i^{\frac{4}{n-2}}\hat{W}_{m,i}+C\hat{W}_{m,i}^{\frac{n+2}{n-2}},
$$where $C$ is some bounded constant which may vary and doesn't depend on $m$. Using Lemma \ref{L5}, Lemma \ref{L1} and Lemma \ref{L2}, integrating in $\Omega_j$ and combining them together, we deduce that for $n \geq 5$
$$
\int_{{\mathbb R}^n}\hat{W}_{m,i}\sigma_i^{\frac{n+2}{n-2}},\quad \int_{{\mathbb R}^n}\hat{W}_{m,i}^{\frac{n+2}{n-2}}\sigma_i\le \frac{C}{(\lambda l)^{n-2}}\le \frac{C}{\lambda^{\frac{n}{2}}}.
$$

\medskip

Similarly, we obtain that for $n\geq 5$
$$
|\l\sigma_i,\sigma_i\r-\int_{{\mathbb R}^n}K_{\lambda}\sigma_i^{\frac{2n}{n-2}}|\le \frac{C}{\lambda^{\beta}}\le \frac{C}{\lambda^{\frac{n}{2}}},
$$
$$
|\sum_{j\neq i}\l\sigma_j,\sigma_i\r|\le \sum_{j\neq i}\frac{C}{|X^j-X^i|^{n-2}}\le \frac{C}{\lambda^{\frac{n}{2}}}.
$$

Therefore we get that for each $i$,
$$
|G_i(0,0)|\le \frac{C}{\lambda^{\frac{n}{2}}}.
$$

Combining the estimates for $F(0,0)$ and $G_i(0,0)$ together, we obtain (\ref{p1.1}).

\medskip

\noindent
{\bf Step 2.} Next, we show that, for $l$ large,
\begin{equation}\label{p1.2}
\|N'(0,0)^{-1}\|\le C.
\end{equation}

\medskip

We consider the equation
\begin{equation}\label{p1.3}
N'(0,0)(\tilde{\phi},\tilde{\epsilon})=(v,\eta).
\end{equation}

For the $v$ component, the equation becomes
\begin{equation}\label{p1.4}
v=\tilde{\phi}-\frac{n+2}{n-2}P_{\mathcal{E}}(-\Delta)^{-1}K_{\lambda}W_m^{\frac{4}{n-2}}\tilde{\phi}-\frac{n+2}{n-2}P_{\mathcal{E}}(-\Delta)^{-1}K_{\lambda}W_m^{\frac{4}{n-2}}(\sum_i\tilde{\epsilon}_i\sigma_i).
\end{equation}

For each $i$-th component, we have
\begin{equation}\label{p1.5}
\eta_i=-\frac{n+2}{n-2}\int_{{\mathbb R}^n}K_{\lambda}W_m^{\frac{4}{n-2}}\sigma_i\tilde{\phi}+\sum_j\tilde{\epsilon_j}\left(\l\sigma_i,\sigma_j\r-\frac{n+2}{n-2}\int_{{\mathbb R}^n}K_{\lambda}W_m^{\frac{4}{n-2}}\sigma_i\sigma_j\right).
\end{equation}

Since $\l v,\sigma_i\r=0$, $\l \tilde{\phi},\sigma_i\r=0$ for all $i$, in (\ref{p1.4}), we can replace $K_{\lambda}W_m^{\frac{4}{n-2}}(\sum_i\tilde{\epsilon}_i\sigma_i)$ by $\tilde{l}_m:=K_{\lambda}W_m^{\frac{4}{n-2}}(\sum_i\tilde{\epsilon_i}\sigma_i)-\sum_i\tilde{\epsilon}_i\sigma_i^{\frac{n+2}{n-2}}$. Then we get
\begin{equation}\label{p1.6}
\begin{array}{ll}
\tilde{\phi}(y)=&v(y)+\frac{n+2}{n(n-2)^2\omega_n}\int_{{\mathbb R}^n}\frac{1}{|z-y|^{n-2}}K_{\lambda}(z)W_m^{\frac{4}{n-2}}(z)\tilde{\phi}(z)dz\\
\\
&+\frac{1}{n(n-2)\omega_n}\int_{{\mathbb R}^n}\frac{\tilde{l}_m(z)}{|y-z|^{n-2}}dz+\sum_{i,j} c_{i,j}Z_{i,j}(y)+\sum_i b_i\sigma_i(y)\\
\\
&=I+II+v+\sum_{i,j} c_{i,j}Z_{i,j}(y)+\sum_i b_i\sigma_i(y).
\end{array}
\end{equation}

\medskip

For $II$, similar to the estimate of $l_m$, we have
$$
\|\tilde{l}_m\|_{**}\le\frac{C|\tilde{\epsilon}|}{\lambda^{\frac{n+2}{2}-\tau}}.
$$

To estimate $c_{i,j}$, we multiply $\sigma_s^{\frac{4}{n-2}}Z_{s,t}$ on both side of (\ref{p1.6}). Modifying the proof of Lemma \ref{La.11}, we infer that
$$
|c_{i,j}|, |b_i|\le\left( C\|\tilde{l}_m\|_{**}+\frac{C}{\lambda^{\frac{n}{2}}}\|\tilde{\phi}\|_*\right)\frac{1}{\lambda^{\tau-1}}.
$$

\medskip

Applying Lemma \ref{La.13}, we obtain that, for $l$ large,
\begin{equation}\label{p1.7}
\|\tilde{\phi}\|_*\le C(\|v\|_*+\frac{|\tilde{\epsilon}|}{\lambda^{\frac{n+2}{2}-\tau}}).
\end{equation}

To estimate (\ref{p1.5}), from Lemma \ref{La.12}, we have
$$
|\int_{{\mathbb R}^n}K_{\lambda}W_m^{\frac{4}{n-2}}\sigma_i\tilde{\phi}|\le \frac{C\|\tilde{\phi}\|_*}{\lambda^{\frac{n-2}{2}+\tau}}.
$$

For $j\neq i$, using Lemma \ref{L1}, we obtain
$$
|\l\sigma_i,\sigma_j\r-\frac{n+2}{n-2}\int_{{\mathbb R}^n}K_{\lambda}W_m^{\frac{4}{n-2}}\sigma_j\sigma_i|\le \frac{C}{|X^i-X^j|^{\frac{n-2}{2}}},
$$and it is easy to get that for $l$ large but independent of $m$,

$$
\l\sigma_i,\sigma_i\r-\frac{n+2}{n-2}\int_{{\mathbb R}^n}W_m^{\frac{4}{n-2}}\sigma_i^2\le -\frac{2}{n-2}\l\sigma_i,\sigma_i\r\le -C.
$$

Therefore, we obtain that
\begin{equation}\label{p1.8}
|\tilde{\epsilon}_i|\le C\left(|\eta|+\frac{\|\tilde{\phi}\|_*}{\lambda^{\frac{n-2}{2}+\tau}}+\frac{|\tilde{\epsilon}|}{(\lambda l)^{\frac{n-2}{2}}}\right).
\end{equation}

Estimates (\ref{p1.7}) and (\ref{p1.8}) yield
\begin{equation}\label{p1.9}
\|(\tilde{\phi},\tilde{\epsilon})\|\le C\left(\|(v,\eta)\|+\frac{\|(\tilde{\phi},\tilde{\epsilon})\|}{\lambda^{\frac{n-2}{2}}}\right),
\end{equation}therefore (\ref{p1.2}) follows when $l$ is chosen large (independent of $m$).

\medskip

\noindent
{\bf Step 3.} Next we estimate $N'(\phi,\epsilon)-N'(0,0)$ when $\|(\phi,\epsilon)\|\le \frac{1}{2}$. We compute the term
$$
(N'(\phi,\epsilon)-N'(0,0))(v,\eta)=(\tilde{v},\tilde{\eta}).
$$

\medskip

For the $\tilde{v}$ component, it holds
\begin{equation}\label{p1.10}
\tilde{v}=-\frac{n+2}{n-2}P_{\mathcal E}(-\Delta)^{-1}\left(K_{\lambda}(v+\sum_i\eta_i\sigma_i)\{\left((\bar{W}_m+\phi)^+\right)^{\frac{4}{n-2}}-W_m^{\frac{4}{n-2}}\}\right),
\end{equation}
 and  for the $i$-th component, we have
\begin{equation}\label{p1.11}
\tilde{\eta}_i=-\frac{n+2}{n-2}\int_{{\mathbb R}^n}K_{\lambda}\sigma_i(v+\sum_j\eta_j\sigma_j)\{\left((\bar{W}_m+\phi)^+\right)^{\frac{4}{n-2}}-W_m^{\frac{4}{n-2}}\}.
\end{equation}

We first consider the term
$$
\|K_{\lambda}v\{\left((\bar{W}_m+\phi)^+\right)^{\frac{4}{n-2}}-W_m^{\frac{4}{n-2}}\}\|_{**}.
$$

\medskip

Set $\epsilon W_m=\sum_i\epsilon_i\sigma_i$ and $\Omega_+:=\{x\in{\mathbb R}^n| u(x)\ge 0\}$ and $\Omega_-=\Omega_+^c$. Define
$$
\chi_{\Omega_-}(x)=\left\{\begin{array}{cc}
1& \mbox{if}\quad x\in\Omega_-\\
\\
0&\mbox{if}\quad x\in\Omega_+.
\end{array}\right.
$$

Then we get
\begin{equation}
 \left((\bar{W}_m+\phi)^+\right)^{\frac{4}{n-2}}-W_m^{\frac{4}{n-2}}=(\bar{W}_m+\phi)^{\frac{4}{n-2}}-W_m^{\frac{4}{n-2}}+\chi_{\Omega_-}|\bar{W}_m+\phi|^{\frac{4}{n-2}}.
\end{equation}

\medskip

Since $\epsilon$ is small, there holds
\begin{equation}\label{p1.11-1}
|(\bar{W}_m+\phi)^{\frac{4}{n-2}}-W_m^{\frac{4}{n-2}}|\le C\left\{\begin{array}{l}
W_m^{\frac{4}{n-2}-1}|\epsilon W_m+\phi|,\quad\mbox{if}\quad W_m\ge |\phi|,\\
\\
|\phi|^{\frac{4}{n-2}},\quad\mbox{if}\quad |\phi|\ge W_m,
\end{array}\right.
\end{equation}

$$
|\phi|^{\frac{4}{n-2}}|v|\le \|\phi\|_*^{\frac{4}{n-2}}\|v\|_*\left(\sum_j\frac{\gamma(y)}{(1+|y-x^j|)^{\frac{n-2}{2}+\tau}}\right)^{\frac{n+2}{n-2}}.
$$

\medskip

Similar to the estimate of $\| l_m \|_{**}$, we use Lemma \ref{L5}. If $y\in B_i\cap\Omega_i$,
$$
\left(\sum_j\frac{\gamma(y)}{(1+|y-x^j|)^{\frac{n-2}{2}+\tau}}\right)^{\frac{n+2}{n-2}}\le C\sum_j\frac{\gamma(y)}{(1+|y-x^j|)^{\frac{n+2}{2}+\tau}}.
$$

If $y\in B_i^c\cap B_{i,m}\cap\Omega_i$, we obtain that
$$\begin{array}{l}
\left(\sum_j\frac{\gamma(y)}{(1+|y-x^j|)^{\frac{n-2}{2}+\tau}}\right)^{\frac{n+2}{n-2}}\le C\sum_j\frac{1}{(1+|y-x^j|)^{\frac{n+2}{2}+\tau+\frac{4}{n-2}(\tau-k)}}\frac{1}{(\lambda l)^{\frac{4}{n-2}k}}\\
\\
\le C\sum_j\frac{1}{(1+|y-x^j|)^{\frac{n+2}{2}+\tau}}\frac{1}{(\lambda l)^{\frac{4}{n-2}k}},
\end{array}
$$as $\tau\ge k$.

If $y\in B_{i,m}^c\cap\Omega_i$,  we have first
$$
\left(\sum_j\frac{1}{(1+|y-x^j|)^{\frac{n-2}{2}+\tau}}\right)^{\frac{n+2}{n-2}}\le C\frac{1}{(1+|y-x^i|)^{\frac{n+2}{2}+\frac{n+2}{n-2}\tau}}m^{\frac{n+2}{n-2}k}.
$$

\medskip

On the other hand
$$
\sum_j\frac{1}{(1+|y-x^j|)^{\frac{n+2}{2}+\tau}}\ge C\frac{1}{(1+|y-x^i|)^{\frac{n+2}{2}+\tau}}(1+C[\frac{m}{2}])^k.
$$

It is not hard to see that when $y\in B_{i,m}^c\cap\Omega_i$, there is a constant $C$ independent of $m$ such that for $\tau\ge k$,
$$
\frac{1}{(1+|y-x^i|)^{\frac{4}{n-2}\tau}}m^{\frac{n+2}{n-2}k}\le C\frac{1}{(\lambda l)^{\frac{4}{n-2}k}}(1+C[\frac{m}{2}])^k
$$which gives

$$
\left(\sum_j\frac{\gamma(y)}{(1+|y-x^j|)^{\frac{n-2}{2}+\tau}}\right)^{\frac{n+2}{n-2}}\le C\sum_j\frac{1}{(1+|y-x^j|)^{\frac{n+2}{2}+\tau}}\frac{1}{(\lambda l)^{\frac{4}{n-2}k}}
$$when  $y\in B_{i,m}^c\cap\Omega_i$. Therefore, we obtain that
$$
\|\chi_{\Omega_-}|\bar{W}_m+\phi|^{\frac{4}{n-2}}v\|_{**},\quad \|\phi^{\frac{4}{n-2}}v\|_{**}\le C\|\phi\|_*^{\frac{4}{n-2}}\|v\|_*.
$$

\medskip

When $W_m>|\phi|$ and for $n\ge 5$,
$$
\begin{array}{l}
\|W_m^{\frac{4}{n-2}-1}|\epsilon W_m+\phi|v\|_{**}\le \|W_m^{\frac{2}{n-2}}\phi^{\frac{2}{n-2}}v\|_{**}+|\epsilon|\|W_m^{\frac{4}{n-2}}v\|_{**}\\
\\
\le C(\|v\|_*\|\phi\|_*^{\frac{2}{n-2}}+|\epsilon|\|v\|_*).
\end{array}
$$

\medskip

Combining the above estimates, we  deduce that
$$
\|K_{\lambda}v\{\left((\bar{W}_m+\phi)^+\right)^{\frac{4}{n-2}}-W_m^{\frac{4}{n-2}}\}\|_{**}\le C\|v\|_*(\|\phi\|_*^{\frac{2}{n-2}}+\|\phi\|_*^{\frac{4}{n-2}}+|\epsilon|).
$$

\medskip

Similarly, we have
$$
\|K_{\lambda}(\sum_j\eta_j\sigma_j)\{\left((\bar{W}_m+\phi)^+\right)^{\frac{4}{n-2}}-W_m^{\frac{4}{n-2}}\}\|_{**}\le C|\eta|(|\epsilon|+\|\phi\|_*+|\epsilon|^{\frac{4}{n-2}}+\|\phi\|_*^{\frac{4}{n-2}}).
$$

\medskip

Hence by Lemma \ref{La.11},  we have
\begin{equation}\label{p1.12}
\begin{array}{l}
\|\tilde{v}\|_*\le C\|v\|_*(\|\phi\|_*^{\frac{2}{n-2}}+|\epsilon|)+C|\eta|(|\epsilon|+\|\phi\|_*+|\epsilon|^{\frac{4}{n-2}}+\|\phi\|_*^{\frac{4}{n-2}})\\
\\
\le C\|(v,\eta)\|(\|\phi\|_*^{\frac{2}{n-2}}+|\epsilon|+\frac{\|\phi\|_*}{\lambda^{\tau-1}}+\frac{|\epsilon|^{\frac{4}{n-2}}}{\lambda^{\tau-1}}+\frac{\|\phi\|_*^{\frac{4}{n-2}}}{\lambda^{\tau-1}}).
\end{array}
\end{equation}

\medskip

Next we estimate (\ref{p1.11}). It is not hard to see that from Lemma \ref{L1}, Lemma \ref{L2} and Lemma \ref{L5}, we can deduce that
$$
\begin{array}{l}
|\int_{{\mathbb R}^n}K_{\lambda}\sigma_i(\eta W_m)\{\left((\bar{W}_m+\phi)^+\right)^{\frac{4}{n-2}}-W_m^{\frac{4}{n-2}}\}|\\
\\
\le C|\eta|(|\epsilon|+|\epsilon|^{\frac{4}{n-2}}+\|\phi\|_*^{\frac{4}{n-2}}+\|\phi\|_*),
\end{array}
$$

$$
\left((\bar{W}_m+\phi)^+\right)^{\frac{4}{n-2}}-W_m^{\frac{4}{n-2}}=(\bar{W}_m+\phi)^{\frac{4}{n-2}}-W_m^{\frac{4}{n-2}}+\chi_{\Omega_-}|\bar{W}_m+\phi|^{\frac{4}{n-2}}.
$$

\medskip

Simple computations give
$$
\begin{array}{l}
|\int_{{\mathbb R}^n}K_{\lambda}\sigma_ivW_m^{\frac{4}{n-2}}|\\
\\
\le \int_{{\mathbb R}^n}\frac{C\|v\|_*}{(1+|y-X^i|)^{n-2}}\sum_j\frac{\gamma(y)}{(1+|y-X^j|)^{\frac{n-2}{2}+\tau}}\left(\sum_j\frac{1}{(1+|y-X^j|)^{n-2}}\right)^{\frac{4}{n-2}}\\
\\
\le \frac{C\|v\|_*}{\lambda^{\tau-1}}\sum_j\int_{\Omega_j\cap B_j}\frac{1}{(1+|y-X^i|)^{n-2}}\frac{1}{(1+|y-X^j|)^{\frac{n}{2}+4}}\\
\\
+\frac{C\|v\|_*}{(\lambda l)^{\frac{4k}{n-2}}}\int_{\cup_s(\Omega_s\cap B_s^c)}\frac{1}{(1+|y-X^i|)^{n-2}}\sum_j\frac{1}{(1+|y-X^j|)^{\frac{n-2}{2}+4+\tau-\frac{4k}{n-2}}}\\
\\
\le \frac{C\|v\|_*}{\lambda^{\tau-1}}.
\end{array}
$$

Similarly we obtain
$$
|\int_{{\mathbb R}^n}K_{\lambda}\sigma_ivW_m^{\frac{4}{n-2}-1}\phi|\le \frac{C\|v\|_*\|\phi\|_*}{\lambda^{2\tau-2}},
$$and
$$
|\int_{{\mathbb R}^n}K_{\lambda}\sigma_iv|\phi|^{\frac{4}{n-2}}|\le  C\frac{\|v\|_*}{\lambda^{\tau-1}}\left(\frac{\|\phi\|_*}{\lambda^{\tau-1}}\right)^{\frac{4}{n-2}}.
$$

\medskip

In view of (\ref{p1.11-1}) and the above estimates, we obtain that

$$
\begin{array}{l}
|\int_{{\mathbb R}^n}K_{\lambda}\sigma_iv\{\left((\bar{W}_m+\phi)^+\right)^{\frac{4}{n-2}}-W_m^{\frac{4}{n-2}}\}|\\
\\
\le \frac{C\|v\|_*}{\lambda^{\tau-1}}(|\epsilon|+\frac{\|\phi\|_*}{\lambda^{\tau-1}}+\left(\frac{\|\phi\|_*}{\lambda^{\tau-1}}\right)^{\frac{4}{n-2}}).
\end{array}
$$

Therefore we have
\begin{equation}\label{p1.13}
\begin{array}{l}
\lambda^{\tau-1}|\tilde{\eta}|\le C\|v\|_*(|\epsilon|+\frac{\|\phi\|_*}{\lambda^{\tau-1}}+\left(\frac{\|\phi\|_*}{\lambda^{\tau-1}}\right)^{\frac{4}{n-2}})\\
\\
+C|\eta|\lambda^{\tau-1}(|\epsilon|+|\epsilon|^{\frac{4}{n-2}}+\|\phi\|_*^{\frac{4}{n-2}}+\|\phi\|_*)\\
\\
\le C\|(v,\eta)\|(|\epsilon|+|\epsilon|^{\frac{4}{n-2}}+\|\phi\|_*^{\frac{4}{n-2}}+\|\phi\|_*).
\end{array}
\end{equation}

Thus by (\ref{p1.12}) and (\ref{p1.13}), we obtain
$$
\|(\tilde{v},\tilde{\eta})\|\le C\|(v,\eta)\|(|\epsilon|+|\epsilon|^{\frac{4}{n-2}}+\|\phi\|_*^{\frac{2}{n-2}}+\|\phi\|_*^{\frac{4}{n-2}}+\|\phi\|_*),
$$
which yields that

\begin{equation}\label{p1.14}
\|N'(\phi,\epsilon)-N'(0,0)\|\le C(|\epsilon|+|\epsilon|^{\frac{4}{n-2}}+\|\phi\|_*^{\frac{2}{n-2}}+\|\phi\|_*^{\frac{4}{n-2}}+\|\phi\|_*).
\end{equation}

\medskip

\noindent
{\bf Step 4.} Set $\theta=\frac{1}{2}$ and $a=\frac{C}{\lambda^{\frac{n+2}{2}-\tau}}$ with $C$ so large that $\|N'(0,0)^{-1}N(0,0)\|\le (1-\theta)a$. For $\|(\phi,\epsilon)\|\le a$, it follows from our estimate (\ref{p1.14}) that
$$
\|N'(0,0)^{-1}\|\|N'(\phi,\epsilon)-N'(0,0)\|\le Ca^{\frac{2}{n-2}}<\theta
$$for $\lambda$ large. The condition of Lemma \ref{L1.1} is satisfied and the existence and uniqueness of $\phi(P,\Lambda)$ and $\epsilon(P,\Lambda)$ follow from the lemma. The $C^1$ dependence follows from (\ref{p1.2}), (\ref{p1.14})and the fact that $N$ has $C^1$ dependence on $P$, $\epsilon$, $\Lambda$ and $\phi$.
\end{proof}

\section{Solving a finite dimensional problem}\label{S2}
\setcounter{equation}{0}

In this section, we will choose the positive constants $C_1$, $C_2$ and integer $l_0$ and solve (\ref{2.9c}) and (\ref{2.9d}) for some $P^i\in B_{\frac{1}{2}}(X^i)$ and $\Lambda_i\in (C_1,C_2)$ when $l\ge l_0$.

To this end, we need some preliminary computations.

\begin{Lemma}\label{L3.a}
\begin{equation}\label{3.1}
(1+\epsilon_i)^{-1}\frac{\partial J}{\partial {P^i}_j}=-\int_{{\mathbb R}^n}K_{\lambda}\sigma_i^{\frac{n+2}{n-2}}Z_{i,j}+o(\frac{1}{\lambda^{\beta}}).
\end{equation}
\end{Lemma}

\begin{proof} The left hand side of (\ref{3.1}) equals
$$
\begin{array}{l}
(1+\epsilon_i)^{-1}\frac{\partial J}{\partial {P^i}_j}=\sum_{s\neq i}(1+\epsilon_s)\l \frac{\partial\sigma_i}{\partial {P^i}_j},\sigma_s\r\\
\\
-\int_{{\mathbb R}^n}K_{\lambda}\left((W_m+\epsilon W_m+\phi)^+\right)^{\frac{n+2}{n-2}}Z_{i,j} \\
\\
=\sum_{s\neq i}(1+\epsilon_s)\l\frac{\partial\sigma_i}{\partial {P^i}_j},\sigma_s\r-\int_{{\mathbb R}^n}K_{\lambda}\left(W_m+\epsilon W_m+\phi\right)^{\frac{n+2}{n-2}}Z_{i,j}\\
\\
-\int_{\Omega_-}K_{\lambda}|W_m+\epsilon W_m+\phi|^{\frac{n+2}{n-2}}Z_{i,j}.
\end{array}
$$

BY Proposition \ref{p2}, $|\epsilon|$ is small. Therefore in $\Omega_-$, $|\phi|\ge \frac{W_m}{2}$.  From Lemma \ref{L8} and the estimates in Lemma \ref{La.2}, we deduce that
$$
\begin{array}{c}
|\int_{\Omega_-}K_{\lambda}|W_m+\epsilon W_m+\phi|^{\frac{n+2}{n-2}}Z_{i,j}|\le C\int_{|\phi|\ge \frac{W_m}{2}}|\phi|^{\frac{n+2}{n-2}}|Z_{i,j}|\\
\\
\le C\|\phi\|_{*}^{\frac{n+2}{n-2}}\frac{1}{(\lambda l)^{\frac{n}{2}+\tau\frac{n+2}{n-2}}}\le  \frac{C}{\lambda^{n+\frac{2(n+2)}{n-2}}}.
\end{array}
$$

Lemma \ref{L3.a} now follows from  Lemma \ref{La.2} and Proposition \ref{p2}.
\end{proof}

\begin{Lemma}\label{L3.b}
\begin{equation}\label{3.3}
\begin{array}{l}
(1+\epsilon_i)^{-1}\frac{\partial J}{\partial\Lambda_i}=\sum_{j\neq
i}\l\frac{\partial\sigma_i}{\partial\Lambda_i},\sigma_j\r\\
\\
-\int_{{\mathbb R}^n} K_{\lambda}\left(\sum_k\sigma_k^{\frac{n+2}{n-2}}+\frac{n+2}{n-2}\sigma_i^{\frac{4}{n-2}}\sum_{j\neq i}\sigma_j \right)Z_{i,n+1}+o(\frac{1}{\lambda^{\beta}}).
\end{array}
\end{equation}
\end{Lemma}
\begin{proof}
As before the left hand side equals
$$
\begin{array}{l}
(1+\epsilon_i)^{-1}\frac{\partial J}{\partial\Lambda_i}=\sum_{j\neq
i}(1+\epsilon_j)\l\frac{\partial\sigma_i}{\partial\Lambda_i},\sigma_j\r\\
\\
-\int_{{\mathbb R}^n} K_{\lambda}\left((\bar{W}_m+\phi|)^+\right)^{\frac{n+2}{n-2}}Z_{i,n+1}\\
\\
=\sum_{j\neq i}(1+\epsilon_j)\l\frac{\partial\sigma_i}{\partial\Lambda_i},\sigma_j\r-\int_{{\mathbb R}^n} K_{\lambda}|\bar{W}_m+\phi|^{\frac{4}{n-2}}(\bar{W}_m+\phi)Z_{i,n+1}\\
\\
-\int_{\Omega_-} K_{\lambda}|\bar{W}_m+\phi|^{\frac{n+2}{n-2}}Z_{i,n+1},
\end{array}
$$
where
$$
\begin{array}{l}
\int_{{\mathbb R}^n} K_{\lambda}|\bar{W}_m+\phi|^{\frac{4}{n-2}}(\bar{W}_m+\phi)Z_{i,n+1}=\int_{{\mathbb R}^n} K_{\lambda}W_m^{\frac{n+2}{n-2}}Z_{i,n+1}\\
\\
+\frac{n+2}{n-2}\int_{{\mathbb R}^n}K_{\lambda}W_m^{\frac{4}{n-2}}(\epsilon W_m+\phi)Z_{i,n+1}\\
\\
+O(1)\int_{|\phi|\ge W_m}|\phi|^{\frac{n+2}{n-2}}|Z_{i,n+1}|+O(1)\int_{|\phi|\le W_m}W_m^{\frac{6-n}{n-2}}|\epsilon W_m+\phi|^2|Z_{i,n+1}|.
\end{array}
$$

Here $O(1)$ means a number uniformly bounded independent of $m$ or $l$.

\medskip

Recall that $|Z_{i,n+1}|\le \frac{C}{(1+|y-P^i|)^{n-2}}$. When $\lambda$ is large, by Lemma \ref{L8}, if $|\phi|\ge W_m$, then  $y\in \left(\cup_j (B_j\cap \Omega_j)\right)^c:=\Omega$. By Lemma \ref{L1} and similar argument as in Lemma \ref{La.2}, we deduce that

$$
|\int_{\Omega_-} K_{\lambda}|\bar{W}_m+\phi|^{\frac{n+2}{n-2}}Z_{i,n+1}|\le C\int_{\Omega}|\phi|^{\frac{n+2}{n-2}}|Z_{i,n+1}|\le \frac{C}{\lambda^{n-1+\frac{2(n+2)}{n-2}}}.
$$

$$
\int_{|\phi|\le W_m}W_m^{\frac{6-n}{n-2}}|\epsilon W_m+\phi|^2 |Z_{i,n+1}|\le C(\frac{\|\phi\|_*^2}{\lambda^{2\tau-2}}+|\epsilon|^2)\le \frac{C}{\lambda^n}.
$$

Since $\l\phi,Z_{i,n+1}\r=0$, by Lemma \ref{La.12},
$$
\int_{{\mathbb R}^n}K_{\lambda}W_m^{\frac{4}{n-2}}\phi Z_{i,n+1}=\int_{{\mathbb R}^n}(K_{\lambda}-1)\sigma_i^{\frac{4}{n-2}}\phi Z_{i,n+1}
+o(\frac{1}{\lambda^n})\le \frac{C}{\lambda^n}.
$$

Similarly
$$
\begin{array}{l}
\int_{{\mathbb R}^n}K_{\lambda}W_m^{\frac{4}{n-2}}\epsilon W_m Z_{i,n+1}\\
\\
=\int_{{\mathbb R}^n}(K_{\lambda}-1)\epsilon_i\sigma_i^{\frac{n+2}{n-2}}Z_{i,n+1}+\frac{C|\epsilon|}{(\lambda l)^{n-2}}\le o(\frac{1}{\lambda^n}).
\end{array}
$$

In the region $y\in \Omega_j\cap B_j^c$
$$
W_m^{\frac{n+2}{n-2}}(y)\le \frac{C}{(\lambda l)^{\frac{4}{n-2}k}}\sum_s\frac{1}{(1+|y-X^s|)^{n+2-\frac{4}{n-2}k}}.
$$
Therefore by Lemma \ref{L1} we have that for all $n\ge 5$
$$\begin{array}{l}
|\int_{\cup_j(\Omega_j\cap B_j^c)}K_{\lambda}(y)W_m^{\frac{n+2}{n-2}}Z_{i,n+1}dy|\le\\
\\
\frac{C}{(\lambda l)^{\frac{4}{n-2}k}}\int_{\cup_j(\Omega_j\cap B_j^c)}\sum_s\frac{1}{(1+|y-X^s|)^{n+2-\frac{4}{n-2}k}}\frac{1}{(1+|y-X^i|)^{n-2}}dy\\
\\
\le \frac{C}{(\lambda l)^{\frac{4}{n-2}k}}\sum_{s\neq  i}\frac{1}{|X^i-X^s|^{\frac{n-2}{2}}}\frac{1}{(\lambda l)^{\frac{n+2}{2}-\frac{4}{n-2}k}}+\frac{C}{(\lambda l)^{n}}\\
\\
\le \frac{C}{(\lambda l)^n}
\end{array}
$$

For $y\in \Omega_j\cap B_j$, we have
$$\begin{array}{l}
|W_m^{\frac{n+2}{n-2}}-\sigma_j^{\frac{n+2}{n-2}}-\frac{n+2}{n-2}\sigma_j^{\frac{4}{n-2}}\hat{W}_{m,j}|\\
\\
\le C \hat{W}_{m,j}^2\sigma_j^{\frac{4}{n-2}-1}\le \hat{W}_{m,j}^{\frac{n}{n-2}}\sigma_j^{\frac{n}{n-2}-1}.
\end{array}
$$

\medskip

We can show for $j\neq i$,
$$\begin{array}{l}
\int_{\Omega_j\cap B_j}K_{\lambda}(y)W_m^{\frac{n+2}{n-2}}Z_{i,n+1}dy=\int_{\Omega_j\cap B_j}K_{\lambda}(y)\sigma_j^{\frac{n+2}{n-2}}Z_{i,n+1}dy\\
\\
+\frac{C}{(\lambda l)^2|X^i-X^j|^{n-2}}=\int_{{\mathbb R}^n}K_{\lambda}(y)\sigma_j^{\frac{n+2}{n-2}}Z_{i,n+1}dy+\frac{C}{(\lambda l)^2|X^i-X^j|^{n-3}}.
\end{array}
$$

When $j=i$,
$$\begin{array}{l}
\int_{\Omega_i\cap B_i}K_{\lambda}(y)W_m^{\frac{n+2}{n-2}}Z_{i,n+1}dy\\
\\
=\int_{\Omega_i\cap B_i}K_{\lambda}(y)\left(\sigma_i^{\frac{n+2}{n-2}}+\frac{n+2}{n-2}\sigma_i^{\frac{4}{n-2}}\sum_{j\neq i}\sigma_j\right)Z_{i,j}dy+\frac{C}{(\lambda l)^{n-1}}\\
\\
=\int_{{\mathbb R}^n}K_{\lambda}(y)\left(\sigma_i^{\frac{n+2}{n-2}}+\frac{n+2}{n-2}\sigma_i^{\frac{4}{n-2}}\sum_{j\neq i}\sigma_j\right)Z_{i,n+1}dy+\frac{C}{(\lambda l)^{n-1}}.
\end{array}
$$

Together with
$$
| \sum_{j\neq i}\epsilon_j \l\frac{\p\sigma_i}{\p \Lambda_i},\sigma_j\r\le \frac{C|\epsilon|}{(\lambda l)^{n-2}} | =o(\frac{1}{\lambda^{\beta}}),
$$
we can easily deduce the estimate (\ref{3.3}).

\end{proof}

\medskip

By the estimates in \cite{B},
$$
\int_{{\mathbb R}^n}\sigma_j^{\frac{n+2}{n-2}}\frac{\partial\sigma_i}{\partial\lambda_i}=C_4\frac{\partial\epsilon_{ij}}{\partial\lambda_i}+\frac{1}{\lambda_i}O(\epsilon_{ij}^{\frac{n}{n-2}}\log\epsilon_{ij}^{-1})
$$where $C_4=\left(n(n-2)\right)^{\frac{n}{2}}\int_{{\mathbb R}^n}\frac{1}{(1+|y|^2)^{\frac{n+2}{2}}}dy$ and

$$
\epsilon_{ij}=\left(\frac{\lambda_i}{\lambda_j}+\frac{\lambda_j}{\lambda_i}+\lambda_i\lambda_j|P^i-P^j|^2\right)^{-\frac{n-2}{2}},\quad\mbox{for $i\neq j$}.
$$

\medskip

Using Lemma \ref{L3.a} and Lemma \ref{La.1}, we infer that (\ref{2.9c}) is equivalent  to
\begin{equation}\label{3.9a}
D_{n,\beta}a_j\frac{1}{\Lambda_i^{\beta-2}\lambda^{\beta}}({P^i}_j-{X^i}_j)=O(\frac{|P^i-X^i|^2}{\lambda^{\beta}})+o(\frac{1}{\lambda^{\beta}}),
\end{equation}for all $i=1,...,(m+1)^k$ and $j=1,...,n$.

\medskip

By Lemma \ref{L3.b}, Lemma \ref{La.3}, Lemma \ref{La.4} and Lemma \ref{La.5}, we can derive that (\ref{2.9d}) is equivalent to

\begin{equation}\label{3.9}
\sum_{j\neq i}\frac{(n-2)C_4A_{ij}\Lambda_j}{2(\Lambda_i\Lambda_j)^{\frac{n}{2}}(\lambda l)^{n-2}}-\frac{C_3}{\Lambda_i^{\beta+1}\lambda^{\beta}}=o(\frac{1}{\lambda^{\beta}})+O(\frac{|P^i-X^i|^{\beta-1}}{\lambda^{\beta}}).
\end{equation}

\medskip

In the above, $A:=\{A_{ij}\}$ is a $(m+1)^k\times (m+1)^k$ matrix associated to the lattice $X_{l,m}$(or $X_{1,m}$), given as follows
$$
A=(A_{ij})=\left\{\begin{array}{l} 0\quad\mbox{if}\quad i=j\\
\\
\left(\frac{\lambda l}{|X^i-X^j|}\right)^{n-2}\quad\mbox{if}\quad i\neq j.
\end{array}\right.
$$

\medskip

If we take $b_i=\Lambda_i^{-\frac{n-2}{2}}$, then we see that (\ref{2.9d}) is equivalent  to
\begin{equation}\label{3.8}
\sum_{j\neq i}A_{ij}b_j-\left(\frac{2C_3}{(n-2)C_4}+o(1)+O(\frac{|P^i-X^i|^{\beta-1}}{\lambda})\right)b_i^{\frac{2\beta}{n-2}-1}=0.
\end{equation}

\medskip

Now we consider the functional $F: {\mathbb R}_+^{(m+1)^k}\to {\mathbb R}$ defined by
\begin{equation}\label{3.10}
F(b)=\frac{1}{2}b^tAb-\frac{C_3}{\beta C_4}\sum_i b_i^{\frac{2\beta}{n-2}}, \ \ \mbox{for} \ b=(b_1,...,b_{(m+1)^k}).
\end{equation}

Since $C_3>0$ and $\beta>n-2$,  the maximum of $F$ will give a solution to the system
$$
F_i(b)=\sum_{j\neq i}A_{ij}b_j-\frac{2C_3}{(n-2)C_4}b_i^{\frac{2\beta}{n-2}-1}=0, \ \ i=1,...,(m+1)^k.
$$

\medskip

Let $\bar{B}_m=(\bar b_1,...,\bar b_{(m+1)^k})$ be a solution to the above system.

\begin{Lemma}\label{L10}
There exist positive constants $C_5<C_6$ independent of $m$,  such that  $C_5\le |\bar b_i|\le C_6$ for all $i=1,...,(m+1)^k$.
\end{Lemma}
\begin{proof}
For each $F_i$, for any integer $m\ge 1$, without loss of generality, we can assume that  $\bar b_1\le \bar b_i\le \bar b_2$ for all $i=1,...,(m+1)^k$. From the equation $F_2(\bar B_m)=0$, summing in $j$, we can get
$$
\bar b_2^{\frac{2\beta}{n-2}-1}\le C\max_{j\neq i}\bar b_j\le C\bar b_2,
$$where  $C=\frac{(n-2)C_4}{2C_3}\sum_{j\neq i}A_{ij}$ and $\sum_{j\neq i}A_{ij}$ can be controlled by $\int_{{\mathbb R}^k}\frac{dy}{1+|y|^{n-2}}$.

\medskip

Similarly from the equation $F_1(\bar B_m)=0$, summing in $j$, we deduce that
$$
\bar b_1^{\frac{2\beta}{n-2}-1}\ge \frac{(n-2)C_4}{2C_3}\sum_{j\neq i}A_{ij} \bar b_1\ge  \frac{(n-2)C_4}{2C_3}\bar b_1.
$$

\medskip

From the above two inequalities, we conclude that $C_5\le \bar b_i\le C_6$ for all $i=1,...,(m+1)^k$.
\end{proof}

\medskip

By the form of $F$ and using the fact that  $2<\frac{2\beta}{n-2}<4$, we see that the Hessian matrix at $\bar B_m$ $D^2F(\bar B_m)$ is negative definite (this ensures that $\bar B_m$ is unique). We will show that  the inverse matrix  of $D^2F(\bar B_m)$ is uniformly bounded independent of $m$ when $l$  is large enough.

Take $X=(x_1,...,x_{(m+1)^k})$ in ${\mathbb R}^{(m+1)^k}$,
$$
\left(D^2(F(\bar B_m))X\right)_i=\sum_{j\neq i}A_{ij}x_j-(\frac{2\beta}{n-2}-1)\frac{2C_3}{(n-2)C_4}\bar b_i^{\frac{2\beta}{n-2}-2}x_i.
$$

Consider the $i$ with largest $|\frac{x_i}{\bar b_i}|$ (from Lemma \ref{L10}, $|x_i|\ge C|X|$). By the equation $F_i(\bar B_m)=0$, as $2\beta>2(n-2)$, we obtain
$$
\begin{array}{l}
|(\frac{2\beta}{n-2}-1)\frac{2C_3}{(n-2)C_4}\bar b_i^{\frac{2\beta}{n-2}-2}x_i|\ge |\frac{2C_3}{(n-2)C_4}\bar b_i^{\frac{2\beta}{n-2}-1}\frac{x_i}{\bar b_i}|\\
\\
\ge|\sum_{j\neq i}A_{ij}\bar b_j|\times|\frac{x_i}{\bar b_i}|\ge \sum_{j\neq i}A_{ij}|x_j|\ge |\sum_{j\neq i}A_{ij}x_j|.
\end{array}
$$

\medskip

This implies that
$$\begin{array}{l}
|\left(D^2(F(\bar B_m))X\right)_i|\ge |(\frac{2\beta}{n-2}-1)\frac{2C_3}{(n-2)C_4}\bar b_i^{\frac{2\beta}{n-2}-2}x_i|-|\sum_{j\neq i}A_{ij}x_j|\\
\\
\ge |\frac{x_i}{\bar b_i}|(\frac{2\beta}{n-2}-1)\frac{2C_3}{(n-2)C_4}\bar b_i^{\frac{2\beta}{n-2}-2}\bar b_i-|\frac{x_i}{\bar b_i}|(\sum_{j\neq i}A_{ij}\bar b_j),\\
\\
=(\frac{2\beta}{n-2}-2)\frac{2C_3}{(n-2)C_4}\bar b_i^{\frac{2\beta}{n-2}-1}|\frac{x_i}{\bar b_i}| \quad (\mbox{by}\quad F_i(\bar B_m)=0)\\
\\
\ge C|X| \quad (\mbox{by Lemma \ref{L10}}),
\end{array}
$$where $C$ only depends on $C_3, C_4, C_5, C_6$.  Hence we get

$$
|D^2F(\bar B_m)X|\ge C|X|.
$$

Similarly, we can also show that
$$
|D^2 F(\bar B_m)X|\le C|X|.
$$

\medskip

Thus we obtain that $|\left(D^2F(\bar B_m)\right)^{-1}X|\le C|X|$ for all $X\in{\mathbb R}^{(m+1)^k}$ with maximum norm $ | \cdot |$.

\begin{proof}[\bf Proof of Theorem \ref{T2}.]
By (\ref{3.9a}), (\ref{2.9c}) is equivalent to
\begin{equation}\label{3.9b}
P^i-X^i=O(|P^i-X^i|^2)+o(1),\quad \mbox{for all}\quad i,
\end{equation}

From (\ref{3.8}) if we let $t=b-\bar B_m\in {\mathbb R}^{(m+1)^k}$, then (\ref{2.9d}) is equivalent to
\begin{equation}\label{3.9c}
D^2F(\bar B_m)t=O(|t|^2)+o(1)+O(\max |P^i-X^i|^2).
\end{equation}

\medskip

For ${\mathbb R}^{(m+1)^k\times (n+1)}$, equipped with maximum norm, we can choose a $C>0$ large but independent of $m$ and $l$.  When $l$ is large enough, (\ref{3.9b}) and (\ref{3.9c}) define a continuous map from
$$
B:=B_{Co(1)}(X^1)\times...\times B_{Co(1)}(X^{(m+1)^k})\times B_{Co(1)}(\bar B_m)\to B.
$$

By Brouwer fixed point theorem, we can solve equations (\ref{3.9b}) and (\ref{3.9c}) near $(X^1,...,X^{(m+1)^k},B_m)$ with
$$
|P^i-X^i|=o(1),\quad |b-B_m|=o(1).
$$

Therefore we have  solved $\frac{\partial J}{\partial \Lambda_i}=0$, $\frac{\partial J}{\partial {P^i}_j}=0$ with
\begin{equation}\label{3.9d}
|\Lambda-B_m^{-\frac{2}{n-2}}|=o(1),\quad |P^i-X^i|=o(1),
\end{equation}
when $\lambda$ large enough.

\medskip

By Lemma \ref{L10}, we can now choose positive constants $C_1<C_2$ which only depend on $C_5$ and $C_6$ and are independent of $m$ and $l$. Then we can take integer $l_0$ large enough such that the $(P,\Lambda)$ given in (\ref{3.9d}) satisfies $P^i\in B_{\frac{1}{2}}(X^i)$ and $\Lambda_i\in (C_1,C_2)$ for all $i$ and $l\ge l_0$.  Therefore a solution to equation (\ref{1.2}) is guaranteed.
\end{proof}

Now we are ready to prove Theorem \ref{T1}.

\medskip

\begin{proof}[\bf Proof of Theorem \ref{T1}.]
Let $\{u_m\}$ denote the solutions of (\ref{1.1} given by Theorem \ref{T2} with $l\ge l_0$ large and fixed. For each $m$, we can find $x_m\in X_{l,m}$, such that
$$
\cup_{m=1}^{\infty}(X_{l,m}-x_m)=X_l^i.
$$

Let
$$
(S_{\lambda}\hat{u}_m)(x)=(S_{\lambda}u_m)(x+x_m).
$$

Then $\hat{u}_m$ satisfies the same equation as $u_m$ due to the periodicity of $K$. We will show that there exists some constant
$C(l)$, independent of $m$, such that
\begin{equation}\label{3.11}
\hat{u}_m(x)\le C(l), \quad \forall\ x\in {\mathbb R}^n.
\end{equation}

Once we have (\ref{3.11}), we then deduce, by elliptic estimates, that
for any $R>1$, there exists some constant
$C_2(l)$, independent of $m$, such that
\begin{equation}\label{3.12}
\|\hat{u}_m\|_{ C^3(B_R) }\le C_2(l), \quad \forall\ m=1,2,3,...
\end{equation}

This implies that we can pass to a subsequence $\{\hat{u}_{m_i}\}$
such that
$$
\hat{u}_{m_i} \to u\ \mbox{in}\ C^2_{loc}( {\mathbb R}^n)
$$
for some non-negative function $u\in C^2( {\mathbb R}^n)$.
Clearly $u$ satisfies
$$
-\Delta u=K(x) u_{+}^{ \frac{n+2}{n-2} }, \quad     \mbox{in}\ {\mathbb R}^n.
$$

By the form of $u$ given by Theorem \ref{T2}, $u$ can not be identically zero, provided that $l$ is large (but independent of $m$). In fact from the form of $u$ given by Theorem \ref{T2}, $u$ is clearly bounded from below in $B_1(0)$ by a positive constant independent of $m$. Therefore, by strong maximum principle, $u>0$ in ${\mathbb R}^n$.

\medskip

It remains to prove  (\ref{3.11}).  This follows from the form of $S_{\lambda}u_m$ (given by Theorem \ref{T2}). In fact by the estimates on $\|\phi\|_*$, when $k<\frac{n-2}{2}$,
$$
\|\phi\|_{ L^\infty({\mathbb R}^n) }\le \|\phi\|_*\sum_i\frac{1}{(1+|x-X^i|)^{\frac{n-2}{2}+\tau}}\le C\|\phi\|_*.
$$where $C$ doesn't depend on $m$. For the same reason, we also have
$$
\|\sum_{X^i\in X_{l,m}} \sigma_{P^i, \Lambda_i}\|_{ L^\infty({\mathbb R}^n)}\le C,
$$
where  $C$ doesn't depend on $m$.

\medskip

Thus it follows from the form of solution given by Theorem \ref{T2}, that
$$
\|S_{\lambda}u_m\|_{ L^\infty({\mathbb R}^n) }\le C.
$$

By the form of $S_{\lambda}$, we get $\|u_m\|_{L^{\infty}({\mathbb R}^n)}\le C\lambda^{\frac{n-2}{2}}=Cl^{\frac{(n-2)^2}{2(\beta-n+2)}}$, (\ref{3.11}) is thus  established.
\end{proof}

\section{Proof of Theorem \ref{T4}}\label{S3}
\setcounter{equation}{0}
We first give a lemma which is used in the proof of Theorem \ref{T4}.
\begin{Lemma}\label{L4-2}
For $n\ge 3$, $0<\alpha<1$, let $f\in C^{\alpha}_{loc}({\mathbb R}^n)$ be nonnegative outside a compact set of ${\mathbb R}^n$. Assume that $u\in C^2({\mathbb R}^n)$ satisfies

$$
-\Delta u=f\quad\mbox{in}\quad {\mathbb R}^n,
$$and
$$
\liminf_{|x|\to\infty}u(x)>-\infty.
$$
Then, for some constant $a\ge \min(0,\liminf_{|x|\to\infty}u(x))$,
$$
u(x)=\frac{1}{n(n-2)\omega_n}\int_{{\mathbb R}^n}\frac{f(\tilde x)d\tilde x}{|x-\tilde x|^{n-2}}+a,\quad \forall x\in{\mathbb R}^n,
$$where $\omega_n$ is the volume of a unit ball in ${\mathbb R}^n$.
\end{Lemma}
\begin{proof}
By adding $-\min(0,\liminf_{|x|\to\infty}u(x))$ to $u$, we may assume, without loss of generality, that $\liminf_{|x|\to\infty}u(x)\ge 0$.

Let
$$
u_i(x)=\frac{1}{n(n-2)\omega_n}\int_{B_i}\frac{f(\tilde x)d\tilde{x}}{|x-\tilde x|^{n-2}},\quad i=1,2,3,....
$$

We know that
\begin{equation}\label{6.6}
\Delta (u-u_i)=0\quad\mbox{in}\quad B_i,
\end{equation}
and, using the fact that $f$ is nonnegative outside a compact set,
$$
u_i\le u_{i+1},\quad \Delta (u-u_i)\le 0,\quad\mbox{in ${\mathbb R}^n$ for large $i$}.
$$

Clearly
$$
\liminf_{|x|\to\infty}(u-u_i)(x)\ge 0.
$$

Thus, by the maximum principle,
$$
u-u_i\ge 0\quad\mbox{in}\quad {\mathbb R}^n\quad\mbox{for large $i$}.
$$

Using the Fatou's Lemma, we obtain
$$
\frac{1}{n(n-2)\omega_n}\int_{{\mathbb R}^n}\frac{f(\tilde x)d\tilde{x}}{|x-\tilde x|^{n-2}}\le \lim_{i\to\infty}u_i(x)\le u(x).
$$

Now, by the Lebesgue  dominated convergence theorem,
$$
\lim_{i\to\infty}u_i(x)=\frac{1}{n(n-2)\omega_n}\int_{{\mathbb R}^n}\frac{f(\tilde x)d\tilde{x}}{|x-\tilde x|^{n-2}}\le u(x),\quad \forall x\in{\mathbb R}^n.
$$

For every $R>0$,
$$
\Delta (u-u_i)=0\quad\mbox{in}\quad B_{2R},\quad\forall i>2R.
$$
We know that $\{u-u_i\}$, for large $i$, is a non-increasing sequence of nonnegative harmonic functions in $B_{2R}$. In particular, $\{u-u_i\}$ is uniformly bounded in $B_{2R}$. By the interior derivative estimates of harmonic functions, the convergence of $\{u-u_i\}$ is $C^2$ in $B_{2R}$. Thus $\{u-u_i\}$ converges to some function $\xi$ in $C^2_{loc}({\mathbb R}^n)$. The entire nonnegative harmonic function $\xi$ is a constant, denoted by $a$. Lemma \ref{L4-2} is established.
\end{proof}

\begin{proof}[\bf Proof of Theorem \ref{T4}]
We prove it by contradiction. Let $u$ be a $C^2$ solution of (\ref{1.1}) satisfying (\ref{A2}) for some $R$, $\epsilon>0$ and $0\le i\le k$. We divide the proof into three steps.

\medskip

\noindent
{\bf Step 1.} For any $a>0$, we have
$$
\sup\{u(x)|x\in {\mathbb R}^n,\quad dist(x,{\mathbb R}_i^k)<a\}<\infty.
$$

Suppose not, by making a translation according to the periods of $K$, we may assume that there exists $|x_j|\le a+1$, such that $u_j$, the corresponding translations of $u$, satisfies
$$
\left\{\begin{array}{l}
-\Delta u_j=Ku_j^{\frac{n+2}{n-2}},\quad u_j>0\quad\mbox{in}\quad {\mathbb R}^n,\\
\\
u_j(x_j)\to\infty
\end{array}\right.
$$ and
\begin{equation}\label{4.4}
\inf_{x\in {\mathbb R}^k_i}\sup_{B_R(x)}u_j\ge \epsilon.
\end{equation}

\medskip

By Theorem 1.2 in \cite{Lin}, there exists a positive constant $C(K,a)$ such that
$$
\int_{B_{2a}(0)}|\nabla u_j|^2+u_j^{\frac{2n}{n-2}} \le C(K,a).
$$
\medskip

By Proposition \ref{p3}, $\{u_j\}$, after passing to a subsequence, only has isolated simple blow up points in ${\mathbb R}^n$. Let $S$ be the set of blow up points in ${\mathbb R}^n$. We know that $S\neq \emptyset$.  Proposition 4.2 of \cite{L2}, applied to translations of $\{u_j\}$, shows that there exists a $\delta>0$, such that
$$
\inf_{x,y\in S,x\neq y}|x-y|\ge\delta.
$$

\medskip

Passing to a subsequence and replacing $x_j$ by some nearby points if necessary, we may assume that $x_j\to \bar{x}\in S$ is an isolated simple blow up point. Thus by Proposition 2.3 in \cite{L2}, that
$$
u_j(x_j)u_j\to h\quad C_{loc}^2({\mathbb R}^n\setminus S),
$$
where $h$ is a positive harmonic function on ${\mathbb R}^n\setminus S$ and has a singularity at each point in $S$. By the proof of Theorem 4.2 in \cite{L2}, $S$ can not have more than one point, so $S=\{\bar{x}\}$, and $u_j\to 0$ uniformly on any compact subset of ${\mathbb R}^n\setminus\{\bar{x}\}$. This contradicts (\ref{4.4}).

\medskip

\noindent
{\bf Step 2.} For any $a>0$,
$$
\inf\{ u(x)| x\in {\mathbb R}^n,\quad dist(x,{\mathbb R}^k_i)<a\}>0.
$$

By step 1,
$$
\sup_{dist(x,{\mathbb R}^k_i)<2a}u(x)=C(a)<\infty.
$$

Since
$$
-\Delta u=Ku^{\frac{n+2}{n-2}}=\left( K u^{\frac{4}{n-2}}\right)u,
$$and $|Ku^{\frac{4}{n-2}}|\le (\sup K )C(a)^{\frac{4}{n-2}}$ if $dist(x,{\mathbb R}^k_i)<2a$, We apply the Harnack inequality to obtain
$$
\sup_{B_a(x)}u\le C(a, \sup K)\inf_{B_a(x)}u,\quad\forall x\in {\mathbb R}^k_i.
$$

\medskip

We may assume without loss of generality that $a\ge R$. Then, in view of (\ref{A2}), we have
$$
\sup_{B_a(x)}u\ge \epsilon,\quad \forall x\in {\mathbb R}^k_i.
$$
Step 2 is established.

Taking $a=1$ in step 2, we have, for some positive constant $b$,
\begin{equation}\label{4.4-2}
u(x)\ge b,\quad  \forall x\mbox{ such that } dist(x,{\mathbb R}^k_i)<1.
\end{equation}

\medskip

\noindent
{\bf Step 3.} When $k\ge \frac{n-2}{2}$, (\ref{1.1}) has no $C^2$ solution satisfying (\ref{4.4-2}).

\medskip

By Lemma \ref{L4-2},
\begin{equation}\label{4.4-1}
u(x)=\frac{1}{n(n-2)\omega_n}\int_{{\mathbb R}^n}\frac{K(\tilde x)u^{\frac{n+2}{n-2}}(\tilde x)}{|x-\tilde x|^{n-2}}d\tilde x+a,\quad\forall x\in{\mathbb R}^n,
\end{equation}
where $a\ge 0$. We will show that $a=0$.

Since $u>0$ in ${\mathbb R}^n$ and $\inf K>0$, (\ref{4.4-1}) implies
$$
u(x)\ge \frac{(\inf K)a^{\frac{n+2}{n-2}}}{n(n-2)\omega_n}\int_{{\mathbb R}^n}\frac{d\tilde x}{|x-\tilde x|^{n-2}},
$$therefore $a=0$, since $\int_{{\mathbb R}^n}\frac{d\tilde{x}}{|x-\tilde x|^{n-2}}=\infty$.

From (\ref{4.4-1}) with $a=0$ and (\ref{4.4-2}), we have, for some constant $C>1$,
\begin{equation}\label{6.12}
u(x)\ge \frac{1}{C}\int_{dist(\tilde x,{\mathbb R}^k_i)<1}\frac{d\tilde x}{|x-\tilde x|^{n-2}},\quad\forall x\in{\mathbb R}^n.
\end{equation}
If $k\ge n-2$, the right hand side of the above is $\infty$, which is impossible.

Now we treat the remaining case: $\frac{n-2}{2}\le k<n-2$. For convenience, we write ${\mathbb R}^n={\mathbb R}^k\times{\mathbb R}^{n-k}$. For any $x\in{\mathbb R}^n$, $x=(y,z)\in{\mathbb R}^k\times{\mathbb R}^{n-k}$ and $u(x)=u(y,z)$. We show that, for some constant $C>1$,
\begin{equation}\label{6.12-1}
u(y,z)\ge \frac{1}{C(1+|z|)^{n-2-k}},\quad\forall (y,z)\in{\mathbb R}_i^k\times{\mathbb R}^{n-k}.
\end{equation}
By (\ref{6.12}),
\begin{equation*}
\begin{array}{ll}
u(y,z)&\ge \frac{1}{C}\int_{{\mathbb R}^k_i}\int_{|\tilde{z}|\le 1}\frac{d\tilde y d\tilde z}{|(y,z)-(\tilde y,\tilde z)|^{n-2}},\\
\\
&\ge \frac{1}{C}\int_{{\mathbb R}^k_i}\frac{d\xi}{(1+|\xi|^2)^{\frac{n-2}{2}}}\int_{|\tilde{z}|\le 1}\frac{d\tilde z}{|z-\tilde z|^{n-k-2}},\\
\\
&\ge \frac{1}{2^iC}\int_{{\mathbb R}^k}\frac{d\xi}{(1+|\xi|^2)^{\frac{n-2}{2}}}\int_{|\tilde{z}|\le 1}\frac{d\tilde z}{|z-\tilde z|^{n-k-2}},\\
\\
&\ge \frac{\omega_{n-k}}{2^iC(|z|+1)^{n-2-k}}\int_{{\mathbb R}^k}\frac{d\xi}{(1+|\xi|^2)^{\frac{n-2}{2}}}\ge \frac{1}{C(1+|z|)^{n-2-k}}.
\end{array}
\end{equation*}
\medskip

Here we have made a change of variables $\xi=\frac{\tilde{y}-y}{|\tilde{z}-z|}$ and have used the fact that for every fixed $\tilde{z}-z\neq 0$ and $y\in {\mathbb R}^k_i$, the set
$$
\{\frac{\tilde{y}-y}{|\tilde{z}-z|}|\tilde{y}\in {\mathbb R}^k_i\}\supset {\mathbb R}^k_i.
$$

Let $v$ be the spherical average of $u$ defined by
$$
v(x)=v(|x|)=v(r)=\frac{1}{|\partial B_r(0)|}\int_{\partial B_r(0)}ud S,
$$then by Jensen's inequality,
$$
-\Delta v\ge (\inf K) v^{\frac{n+2}{n-2}},\quad\mbox{in}\quad {\mathbb R}^n.
$$

By some elementary argument, see e. g. \cite{CGS}, we have, for some constant $C>0$,
\begin{equation}\label{6.12-2}
v(x)\le \frac{C}{(1+|x|)^{\frac{n-2}{2}}},\quad\mbox{for any}\quad x\in{\mathbb R}^n.
\end{equation}
For $k>\frac{n-2}{2}$, we obtain ($r=|x|$), using (\ref{6.12-1}) and (\ref{6.12-2}),
\begin{equation*}
\begin{array}{l}
\frac{C}{(1+r)^{\frac{n-2}{2}}}\ge v(r)=\frac{1}{|\partial B_r(0)|}\int_{\partial B_r(0)}ud S\\
\\
\ge\frac{1}{Cr^{n-1}}\int_0^r da\int_{\{|y|=a\}\cap{\mathbb R}^k_i,|z|=\sqrt{r^2-a^2}}u(y,z)\\
\\
\ge \frac{1}{2^iCr^{n-1}}\int_0^r \frac{1}{(1+\sqrt{r^2-a^2})^{n-2-k}}a^{k-1}(\sqrt{r^2-a^2})^{n-k-1}da\\
\\
\ge \frac{1}{C}\int_0^1\frac{s^{k-1}(\sqrt{1-s^2})^{n-k-1}}{(1+r\sqrt{1-s^2})^{n-2-k}} ds\\
\\
\ge \frac{1}{Cr^{n-2-k}}\int_0^1\frac{s^{k-1}(\sqrt{1-s^2})^{n-k-1}}{(1+\sqrt{1-s^2})^{n-2-k}}ds,\quad\mbox{for $r\ge 1$}\\
\\
\ge\frac{1}{C(1+r)^{n-2-k}},\quad\mbox{for $r\ge 1$}.
\end{array}
\end{equation*}
Sending $r$ to $\infty$ leads to a contradiction.

\medskip

For $k=\frac{n-2}{2}$ and $0\le i\le k=\frac{n-2}{2}$, we derive from (\ref{4.4-1}) (notice that $a=0$) and (\ref{6.12-1}) that, for any $(y,z)\in {\mathbb R}_i^k\times {\mathbb R}^{n-k}$,
\begin{equation}\label{6.12-4}
\begin{array}{ll}
u(y,z)\ge \frac{1}{C}\int_{{\mathbb R}_i^k\times{\mathbb R}^{n-k}}\frac{d\tilde{y}d\tilde{z}}{|(y,z)-(\tilde{y},\tilde{z})|^{n-2}(1+|\tilde{z}|)^{(n-2-k)\frac{n+2}{n-2}}}\\
\\
\ge \frac{1}{C}\int_{{\mathbb R}^k_i}\frac{d\xi}{(1+|\xi|^2)^{\frac{n-2}{2}}}\int_{{\mathbb R}^{n-k}}\frac{d\tilde z}{|\tilde z-z|^{n-k-2}(1+|\tilde z|)^{(n-2-k)\frac{n+2}{n-2}}}\\
\\
\ge \frac{\max(1,\log|z|)}{C(1+|z|)^{\frac{n-2}{2}}}\quad\mbox{by Lemma \ref{L2}}.
\end{array}
\end{equation}

By (\ref{6.12-2}) and (\ref{6.12-4}), we obtain when $k=\frac{n-2}{2}$ and $0\le i\le k$,
\begin{equation*}
\begin{array}{l}
\frac{1}{(1+r)^{\frac{n-2}{2}}}\ge v(r)=\frac{1}{|\partial B_r(0)|}\int_{\partial B_r(0)}ud S\\
\\
\ge\frac{1}{Cr^{n-1}}\int_0^r da\int_{\{|y|=a\}\cap{\mathbb R}^k_i,|z|=\sqrt{r^2-a^2}}u(y,z)\\
\\
\ge \frac{1}{2^iCr^{n-1}}\int_0^r \frac{\max(1,\log{r\sqrt{1-a^2}})}{(1+\sqrt{r^2-a^2})^{\frac{n-2}{2}}}a^{\frac{n-2}{2}-1}(\sqrt{r^2-a^2})^{\frac{n}{4}}da\\
\\
\ge \frac{1}{C}\int_0^1\frac{\max(1,\log r\sqrt{1-s^2})}{(1+r\sqrt{1-s^2})^{\frac{n-2}{2}}}s^{\frac{n-2}{2}-1}(\sqrt{1-s^2})^{\frac{n}{4}} ds\\
\\
\ge \frac{1}{C}\int_0^{\frac{1}{2}}\frac{\log \frac{\sqrt3r}{2}}{(1+r\sqrt{1-s^2})^{\frac{n-2}{2}}}s^{\frac{n-2}{2}-1}(\sqrt{1-s^2})^{\frac{n}{4}} ds,\quad\mbox{for $r\ge 10$}\\
\\
\ge\frac{\log \frac{r}{2}}{(1+r)^{\frac{n-2}{2}}},\quad\mbox{for $r\ge 10$}.
\end{array}
\end{equation*}
We also arrive at a contradiction when $r\to\infty$. Thus Theorem \ref{T4} is proved.
\end{proof}

From the proof of Theorem \ref{T4}, it is easy to see that when $k<\frac{n-2}{2}$, $0\le i\le k$ and $K$ has a positive lower bound, (\ref{1.1}) does not admit a solution $u$ satisfying
$$
\lim_{|z|\to\infty}(1+|z|)^{\frac{n-2}{2}}u(y,z)=\infty,\quad\mbox{uniformly in}\quad y\in{\mathbb R}^k_i.
$$

In some sense, $\frac{n-2}{2}$ is a threshold value for the decay power of solutions of (\ref{1.1}).

\begin{Lemma}\label{L4-4}
 Let $1\le k<\frac{n-2}{2}$, suppose that $K\ge 0$, but not identically equal zero and $K$ is bounded from above. Let $u$ be a positive solution of (\ref{1.1}). Assume, for some constants $\tau>0$, that
\begin{equation}\label{6.13}
\sup_{(y,z)\in{\mathbb R}^n}(1+|z|)^{\frac{n-2}{2}+\tau}u(y,z)< \infty.
\end{equation}

Then
\begin{equation}\label{6.14}
\sup_{(y,z)\in{\mathbb R}^n}(1+|z|)^{n-2-k}u(y,z)< \infty.
\end{equation}
\end{Lemma}
\begin{proof}
When $\tau\ge \frac{n-2}{2}-k$, (\ref{6.14}) is obvious. Now we consider the case $0<\tau<\frac{n-2}{2}-k$. By (\ref{4.4-1}) (notice that $a=0$) and (\ref{6.13}), we obtain that for some constant $C>0$,
$$\begin{array}{ll}
u(y,z)&\le C\int_{{\mathbb R}^k}\int_{{\mathbb R}^{n-k}}\frac{1}{|(y,z)-(\tilde y,\tilde z)|^{n-2}}\frac{1}{(1+|\tilde z|)^{\frac{n-2}{2}+\frac{n+2}{n-2}\tau+2}}d\tilde y d\tilde z,\\
\\
&\le C\int_{{\mathbb R}^{n-k}}\frac{1}{|z-\tilde z|^{n-k-2}}\frac{1}{(1+|\tilde z|)^{\frac{n-2}{2}+\frac{n+2}{n-2}\tau+2}}d\tilde z.
\end{array}
$$

Therefore if $\frac{n-2}{2}+\frac{n+2}{n-2}\tau\neq n-k-2$, applying the first part of Lemma \ref{L2}, we have
$$
u(y,z)\le \frac{C}{(1+|z|)^{\min(n-k-2,\frac{n-2}{2}+\frac{n+2}{n-2}\tau)}}.
$$

If $\frac{n-2}{2}+\frac{n+2}{n-2}\tau> n-k-2$, we are done. Therefore we only need to consider the case $\frac{n-2}{2}+\frac{n+2}{n-2}\tau< n-k-2$ if they are not equal. In this case, we get
$$
u(y,z)\le\frac{C}{(1+|z|)^{\frac{n-2}{2}+\frac{n+2}{n-2}\tau}}.
$$

Let $\tau_1=\frac{n+2}{n-2}\tau$ and $\tau_i=\frac{n+2}{n-2}\tau_{i-1}$.  Obviously, $\{\tau_i\}$ is an increasing sequence and hence we can iterate till we get $\frac{n-2}{2}+\frac{n+2}{n-2}\tau_i\ge n-k-2$. If $\frac{n-2}{2}+\frac{n+2}{n-2}\tau_i> n-k-2$, then we are done. If $\frac{n-2}{2}+\frac{n+2}{n-2}\tau_i= n-k-2$, i. e., $\tau_i=(\frac{n-2}{2}-k)\frac{n-2}{n+2}$, we can apply the second part of Lemma \ref{L2} to get
\begin{equation}\label{6.15}
u(y,z)\le C\frac{\max(1,\log|z|)}{(1+|z|)^{n-k-2}},\quad \forall (y,z)\in{\mathbb R}^n.
\end{equation}

Choose any $\tau_{i+1}\in (\tau_i,\frac{n-2}{2}-k)$, then $\frac{n-2}{2}+\frac{n+2}{n-2}\tau_{i+1}\ge n-k-2$. By (\ref{6.15}), we have
$$
u(y,z)\le \frac{C}{(1+|z|)^{\frac{n-2}{2}+\tau_{i+1}}}.
$$
Since $\frac{n-2}{2}+\frac{n+2}{n-2}\tau_{i+1}> n-k-2$, by one more iteration, we get the conclusion.
\end{proof}

\medskip

\begin{remark}
\em
We easily see from the proof of Theorem \ref{T2} and Lemma \ref{4.4} that solutions constructed in Theorem \ref{T2} and Theorem \ref{T1} satisfy (\ref{6.14}).
\end{remark}

\medskip

In the following, for every $k\in[1,\frac{n-2}{2})$, we give examples of positive smooth function $K$ such that  there is a solution $u$ of (\ref{1.1}) satisfying (\ref{A2}) for some $R$, $\epsilon>0$ and all $i\in [0,k]$.

\medskip

Let $(y, z) \in {\mathbb R}^k \times {\mathbb R}^{n-k}$ and $ u(y,z)=v(z)=\frac{1}{(1+|z|^2)^{\frac{n-2}{4}}}$. Direct calculation shows that
$$
-\Delta u(y,z)=-\Delta_z v(z)=\frac{n-2}{2}\left(\frac{n-2}{2}-k+\frac{n-2}{2(1+|z|^2)}\right) u(y,z)^{\frac{n+2}{n-2}}.
$$
Moreover $u$ decays like $\frac{1}{(1+|z|)^{\frac{n-2}{2}}}$. Here $ K(y, z)= \frac{n-2}{2}\left(\frac{n-2}{2}-k+\frac{n-2}{2(1+|z|^2)}\right) $ which is periodic in $y$.

\section{Appendix A}\label{S4}
\setcounter{equation}{0}

In this section, we present the proof of some technical lemmas.

For $x_i$ $x_j$, $y\in {\mathbb R}^n$,  define
$$
g_{ij}(y)=\frac{1}{(1+|y-x_i|)^{\alpha}(1+|y-x_j|)^{\beta}}
$$
where $x_i\neq x_j$ and $\alpha>0$ and $\beta>0$ are two constants.

\medskip

We first prove a lemma which slightly improves the Lemma B.1 in \cite{WY}.

\medskip

\renewcommand{\theLemma}{A.1}
\begin{Lemma}\label{L1}
For any constant $\tau\in [0,\min(\alpha,\beta)]$, we have
$$
g_{ij}(y)\le \frac{2^{\tau}}{(1+|x_i-x_j|)^{\tau}}\left(\frac{1}{(1+|y-x_i|)^{\alpha+\beta-\tau}}+\frac{1}{(1+|y-x_j|)^{\alpha+\beta-\tau}}\right).
$$
\end{Lemma}
\begin{proof}
Let $d=|x_i-x_j|$. If $y\in B_{\frac{d}{2}}(x_i)$, then
$$
|y-x_j|\ge\frac{d}{2},\quad |y-x_j|\ge|y-x_i|,
$$
which implies
$$
g_{ij}(y)\le \frac{1}{(1+\frac{1}{2}d)^{\tau}}\frac{1}{(1+|y-x_i|)^{\alpha+\beta-\tau}},\quad y\in B_{\frac{1}{2}d}(x_i).
$$

Similarly, we have
$$
g_{ij}(y)\le \frac{1}{(1+\frac{1}{2}d)^{\tau}}\frac{1}{(1+|y-x_j|)^{\alpha+\beta-\tau}},\quad y\in B_{\frac{1}{2}d}(x_j).
$$

Now we consider $y\in {\mathbb R}^n\setminus \left(B_{\frac{1}{2}d}(x_i)\cup B_{\frac{1}{2}d}(x_j)\right)$. Then we have $|y-x_i|\ge d$, $|y-x_j|\ge d$. We may also assume that $|y-x_i|\ge |y-x_j|$. This yields that
$$
g_{ij}(y)\le \frac{1}{(1+d)^{\tau}}\frac{1}{(1+|y-x_j|)^{\alpha+\beta-\tau}}.
$$
The result of the Lemma follows easily from above inequalities.
\end{proof}

\bigskip

\renewcommand{\theLemma}{A.2}
\begin{Lemma}\cite{WY}\label{L2}
For any constant $0<\tau$ with $\tau\neq n-2$, there exists a constant $C=C(n,\tau)>1$ such that
$$
\frac{1}{C(1+|y|)^{\min(\tau,n-2)}}\le\int_{{\mathbb R}^n}\frac{1}{|y-z|^{n-2}(1+|z|)^{2+\tau}}dz\le \frac{C}{(1+|y|)^{\min(\tau,n-2)}}.
$$

When $\tau=n-2$, there exists a constant $C=C(n)>1$ such that
$$
\frac{\max(1,\log|y|)}{C(1+|y|)^{n-2}}\le\int_{{\mathbb R}^n}\frac{1}{|y-z|^{n-2}(1+|z|)^n}dz\le \frac{C\max(1,\log|y|)}{(1+|y|)^{n-2}}.
$$
\end{Lemma}
\begin{proof}
This follows from a simple modification of the proof of Lemma B.2 in \cite{WY}. So we omit the details.
\end{proof}

\medskip

Recall that for $X_{l,m}=\{X^i\}_{i=1}^{(m+1)^k}$, $\Omega_i=\{y\in {\mathbb R}^n,\mbox{ such that } |y-X^i|\le |y-X^j|,\mbox{ for all } j\neq i\}$,  $B_i=B_{\lambda l}(X^i)$ and $B_{i,m}=B_{\max(\frac{m}{4},1)\lambda l}(X^i)$.

\medskip

The following lemma provides basic estimates and will be used frequently in the sequel.

\medskip

\renewcommand{\theLemma}{A.3}
\begin{Lemma}\label{L5}
For any $\theta>k$, there exists a constant $C(\theta,k,n)>1$ independent of $m$, such that if $y\in B_i\cap\Omega_i$,

\begin{equation*}\tag{A1}\label{a1-1}
\frac{1}{(1+|y-X^i|)^{\theta}}\le\sum_j\frac{1}{(1+|y-X^j|)^{\theta}}\le \frac{C}{(1+|y-X^i|)^{\theta}};
\end{equation*}

If $y\in B_i^c\cap B_{i,m}\cap\Omega_i$,
\begin{equation*}\tag{A2}\label{a1-2}
\frac{1}{C(1+|y-X^i|)^{\theta-k}(\lambda l)^k}
\le\sum_j\frac{1}{(1+|y-X^j|)^{\theta}}\le \frac{C}{(1+|y-X^i|)^{\theta-k}(\lambda l)^k};
\end{equation*}

and if $y\in B_{i,m}^c\cap\Omega_i$,

\begin{equation*}\tag{A3}\label{a1-3}
\frac{m^k}{C(1+|y-X^i|)^{\theta}}\le\sum_j\frac{1}{(1+|y-X^j|)^{\theta}}\le \frac{Cm^k}{(1+|y-X^i|)^{\theta}}\le \frac{C}{(1+|y-X^i|)^{\theta-k}(\lambda l)^k}.
\end{equation*}
\end{Lemma}

\begin{proof}
For any $y\in \Omega_i$, since $y$ is closest to $X^i$, by triangle inequality, we have
$$
3|y-X^j|\ge |y-X^i|+|X^i-X^j|\ge |y-X^j|.
$$ Hence
\begin{equation*}\label{4.1}\tag{A4}
\begin{array}{c}
\sum_j\frac{1}{(1+|y-X^j|)^{\theta}}\le \frac{C}{(1+|y-X^i|)^{\theta}}\sum_j\frac{1}{\left(1+\frac{|X^i-X^j|}{(1+|y-X^i|)}\right)^{\theta}}\\
\\
\le\frac{C}{(1+|y-X^i|)^{\theta}}\left(1+\int_{[-m-1,m+1]^k}\frac{1}{\left(1+\frac{\lambda l}{1+|y-X^i|}|z|\right)^{\theta}}dz\right)\\
\\
\le\frac{C}{(1+|y-X^i|)^{\theta}}\left(1+ \frac{(1+|y-X^i|)^k}{(\lambda l)^k}\int_{|z|\le\frac{(m+1)\lambda l}{(1+|y-X^i|)}}\frac{1}{(1+|z|)^{\theta}}dz\right)\\
\\
\le \frac{C}{(1+|y-X^i|)^{\theta}}\left(1+ \frac{(1+|y-X^i|)^k}{(\lambda l)^k}\right),\quad \mbox{if}\quad \theta>k.
\end{array}
\end{equation*}

If $y\in B_i\cap\Omega_i$, the inequalities in (\ref{a1-1}) can be easily obtained from the above.

If $y\in B_i^c\cap\Omega_i$, we have

\begin{equation*}\label{4.2}\tag{A5}
\begin{array}{l}
\sum_j\frac{1}{(1+|y-X^j|)^{\theta}}\ge \frac{C}{(1+|y-X^i|)^{\theta}}\sum_j\frac{1}{\left(1+\frac{|X^i-X^j|}{(1+|y-X^i|)}\right)^{\theta}}\\
\\
\ge  \frac{C}{(1+|y-X^i|)^{\theta}}\left(1+
2^{-k}\int_{[0,[\frac{m}{2}]+1]^k\setminus [0,1]^k}\frac{1}{(1+\frac{\lambda l}{(1+|y-X^i|)}|z|)^{\theta}}dz\right)\\
\\
\ge  \frac{C}{(1+|y-X^i|)^{\theta}}\left(1+ \frac{(1+|y-X^i|)^k}{(\lambda l)^k}\int_{\frac{\lambda l}{(1+|y-X^i|)}\le|z|\le \frac{([\frac{m}{2}]+1)\lambda l}{(1+|y-X^i|)}}\frac{1}{(1+|z|)^{\theta}}dz\right),
\end{array}
\end{equation*}
where the constant $2^{-k}$ is due to the reason that for any point $X^j\in X_{l,m}$, the integral region always contains a quadrant of $[0,[\frac{m}{2}]+1]^k\setminus [0,1]^k$ if it is not empty.

Now when $y\in B_i^c\cap B_{i,m}\cap\Omega_i$, we may assume that $m\ge 8$, since $1<[\frac{m}{4}]\le \frac{[\frac{m}{2}]+1}{2}$, we have

$$
\int_{\frac{\lambda l}{(1+|y-X^i|)}\le|z|\le \frac{(\frac{m}{2}]+1)\lambda l}{(1+|y-X^i|)}}\frac{1}{(1+|z|)^{\theta}}dz\ge \int_{1\le |z|\le 2}\frac{1}{(1+|z|)^{\theta}}dz>0.
$$

The inequalities in (\ref{a1-2}) follows easily from above observation and (\ref{4.1}).

When $y\in B_{i,m}^c\cap\Omega_i$, since $\frac{\lambda l}{1+|y-X^i|}\le \frac{4}{m}$,  we have
$$ \int_{[0,[\frac{m}{2}]+1]^k\setminus [0,1]^k}\frac{1}{(1+\frac{\lambda l}{(1+|y-X^i|)}|z|)^{\theta}}dz  \geq  \int_{[0,[\frac{m}{2}]+1]^k\setminus [0,1]^k}\frac{1}{(1+\frac{4}{m} |z|)^{\theta}}dz  \geq C m^k.
$$
Then the inequalities in (\ref{a1-3}) follows from (\ref{4.1}) and (\ref{4.2}).
\end{proof}

\medskip

\renewcommand{\theLemma}{A.4}
\begin{Lemma}\label{L3}
Suppose that $n\ge 5$, $1\le k<\frac{n-2}{2}$ and $0<C_1<C_2<\infty$. We can find a positive constant $\tau_0=\tau_0(n,k)\in(k,\frac{n-2}{2}]$, such that for any $k\le \tau<\tau_0$, there exist constants $\theta=\theta(\tau,n,k)>0$ and $C=C(C_1,C_2,k,n)$, such that
$$
\begin{array}{l}
\int_{{\mathbb R}^n}\frac{1}{|y-z|^{n-2}}W_m^{\frac{4}{n-2}}\gamma(z)\sum_j\frac{1}{(1+|z-X^j|)^{\frac{n-2}{2}+\tau}}dz\le \\
\\
C\gamma(y)\sum_j\frac{1}{(1+|y-X^j|)^{\frac{n-2}{2}+\tau+\theta}}+\frac{C}{(\lambda l)^{\frac{4k}{n-2}}}\gamma(y)\sum_j\frac{1}{(1+|y-X^j|)^{\frac{n-2}{2}+\tau}}.
\end{array}
$$
\end{Lemma}

\medskip

\begin{proof}
Since $n\ge 5$ and $k<\frac{n-2}{2}$, using Lemma \ref{L5}, we obtain that for $z \in B_i \cap \Omega_i$,
$$
W_m^{\frac{4}{n-2}}\sum_j\frac{1}{(1+|z-X^j|)^{\frac{n-2}{2}+\tau}}\le C \frac{1}{(1+|z-X^i|)^{\frac{n+2}{2}+2+\tau}}
$$
and for $z \in B_i^c \cap \Omega_i \cap B_{i, m}$
$$
W_m^{\frac{4}{n-2}}\sum_j\frac{1}{(1+|z-X^j|)^{\frac{n-2}{2}+\tau}}\le C
\frac{1}{(\lambda l)^{\frac{4}{n-2}k}}\sum_j\frac{1}{(1+|z-X^j|)^{\frac{n-2}{2}+4+\tau-\frac{4}{n-2}k}}.
$$
For $z \in  B_{i, m}^c \cap \Omega_i$, we also have
$$ W_m^{\frac{4}{n-2}}\sum_j\frac{1}{(1+|z-X^j|)^{\frac{n-2}{2}+\tau}}
 \leq C \frac{ m^{k+ \frac{4}{n-2} k}}{ (1+ |z-X^i|)^{ \frac{n-2}{2}+ \tau +4 }}
$$
$$
\leq C \frac{1}{(\lambda l)^{\frac{4}{n-2}k}}\sum_j\frac{1}{(1+|z-X^j|)^{\frac{n-2}{2}+4+\tau-\frac{4}{n-2}k}}.
$$

Now we compute
\begin{equation*}
\begin{array}{l}
\int_{\Omega_s\cap B_s}\frac{1}{|y-z|^{n-2}}W_m^{\frac{4}{n-2}}\gamma(z)\sum_j\frac{1}{(1+|z-X^j|)^{\frac{n-2}{2}+\tau}}dz\le\\
\\
 \int_{\Omega_s\cap
B_s}\frac{1}{|y-z|^{n-2}}\left(\frac{1+|z-X^s|}{\lambda}\right)^{\tau-1}\frac{C}{(1+|z-X^s|)^{2+\frac{n+2}{2}+\tau}}dz\\
\\
\le \frac{C\lambda^{1-\tau}}{(1+|y-X^s|)^{\min(\frac{n-2}{2}+2+1,n-2)}}\le \left(\frac{1+|y-X^s|}{\lambda}\right)^{\tau-1}\frac{C}{(1+|y-X^s|)^{\frac{n-2}{2}+\tau+\theta}}\end{array}
\end{equation*}
where $0<\theta<\min(2,\frac{n-2}{2}-1):=\theta_1$.

Similarly we also have
$$
\begin{array}{l}
\int_{\Omega_s\cap B_s}\frac{1}{|y-z|^{n-2}}W_m^{\frac{4}{n-2}}\gamma(z)\sum_j\frac{1}{(1+|z-X^j|)^{\frac{n-2}{2}+\tau}}dz\le\\
\\
\int_{\Omega_s\cap
B_s}\frac{1}{|y-z|^{n-2}}\frac{1}{(1+|z-X^s|)^{2+\frac{n+2}{2}+\tau}}dz\le \frac{C}{(1+|y-X^s|)^{\frac{n-2}{2}+\tau+\theta}},
\end{array}
$$for $0<\theta<\min(2,\frac{n-2}{2}-\tau)):=\theta_2$.

If $y\in \Omega_i\cap B_i$ for some $i$, from the above two inequalities, taking a $\theta\in (0,\min(\theta_1,\theta_2))$, we have
$$
\begin{array}{l}
\int_{\cup_s \Omega_s\cap B_s}\frac{1}{|y-z|^{n-2}}W_m^{\frac{4}{n-2}}\gamma(z)\sum_j\frac{1}{(1+|z-X^j|)^{\frac{n-2}{2}+\tau}}dz\le  \frac{C\gamma(y)}{(1+|y-X^i|)^{\frac{n-2}{2}+\tau+\theta}} \\
\\
+ C\min\left(\frac{1}{\lambda^{\tau-1}}\sum_{s\neq i}\frac{1}{(1+|y-X^s|)^{\frac{n-2}{2}+\theta+1}},\sum_{s\neq i}\frac{1}{(1+|y-X^s|)^{\frac{n-2}{2}+\tau+\theta}}\right)\\
\\
\le \frac{C\gamma(y)}{(1+|y-X^i|)^{\frac{n-2}{2}+\tau+\theta}},\quad\mbox{by Lemma \ref{L5}}\\
\\
\le C\gamma(y)\sum_j\frac{1}{(1+|y-X^j|)^{\frac{n-2}{2}+\tau+\theta}}.
\end{array}
$$

If $y\in \cup_i(\Omega_i\cap B_i^c)$, then $\gamma(y)=1$ and it is easy to see that
$$\begin{array}{l}
\int_{\cup_s (\Omega_s\cap B_s)}\frac{1}{|y-z|^{n-2}}W_m^{\frac{4}{n-2}}\gamma(z)\sum_j\frac{1}{(1+|z-X^j|)^{\frac{n-2}{2}+\tau}}dz\\
\\
\le \frac{C}{(\lambda l)^{\frac{4k}{n-2}}}\int_{{\mathbb R}^n}\frac{1}{|y-z|^{n-2}}\sum_j\frac{1}{(1+|z-X^j|)^{\frac{n-2}{2}+\tau+4-\frac{4k}{n-2}}}\\
\\
\le \sum_j\frac{C}{(1+|y-X^j|)^{\frac{n-2}{2}+\tau+\theta}},
\end{array}
$$
where $0<\theta<\min(2-\frac{4k}{n-2},\frac{n-2}{2}-\tau):=\theta_3$.

\medskip

Thus we get
\begin{equation*}
\begin{array}{l}
\int_{\cup_i(\Omega_i\cap B_i)}\frac{1}{|y-z|^{n-2}}W_m^{\frac{4}{n-2}}\gamma(z)\sum_j\frac{1}{(1+|z-X^j|)^{\frac{n-2}{2}+\tau}}dz\le \\
\\
C\gamma(y)\sum_j\frac{1}{(1+|y-X^j|)^{\frac{n-2}{2}+\tau+\theta}},
\end{array}
\end{equation*} for all $0<\theta<\min(\theta_1,\theta_2,\theta_3)$.

When $z\in\Omega_i \cap B_i^c$, we estimate as follows: if $y\in \cup_i(\Omega_i\cap B_i^c)$, i.e., $\gamma(y)=1$, we have
$$
\begin{array}{l}
\int_{\cup_i(\Omega_i\cap B_i^c)}\frac{1}{|y-z|^{n-2}}W_m^{\frac{4}{n-2}}\gamma(z)\sum_j\frac{1}{(1+|z-X^j|)^{\frac{n-2}{2}+\tau}}dz \\
\\
\le C\frac{1}{(\lambda l)^{\frac{4}{n-2}k}}\int_{{\mathbb R}^n}\frac{1}{|y-z|^{n-2}}\sum_j\frac{1}{(1+|z-X^j|)^{4+\frac{n-2}{2}+\tau-\frac{4}{n-2}k}}dz\\
\\
\le C\sum_j\frac{1}{(1+|y-X^j|)^{\min(\frac{n-2}{2}+2+\tau-\frac{4}{n-2}k,n-2)}}\frac{1}{(\lambda l)^{\frac{4}{n-2}k}}\\
\\
\le C\left\{\begin{array}{l}
\sum_j\frac{1}{(1+|y-X^j|)^{\frac{n-2}{2}+\tau}}\frac{1}{(\lambda l)^{\frac{4k}{n-2}}},\quad\mbox{when}\quad n\ge 6\\
\\
\sum_j\frac{1}{(1+|y-X^j|)^{n-2}}\frac{1}{(\lambda l)^{\frac{4k}{n-2}}},\quad\mbox{when}\quad n=5, k=1,
\end{array}\right.\\
\\
\le C\gamma(y)\sum_j\frac{1}{(1+|y-X^j|)^{\frac{n-2}{2}+\tau}}\frac{1}{(\lambda l)^{\frac{4k}{n-2}}}.
\end{array}
$$

\medskip

Here in the case $n\ge 6$, we need $\frac{n-2}{2}+2+\tau-\frac{4}{n-2}k<n-2$ which gives $\tau<\frac{n-2}{2}-2+\frac{4}{n-2}k=\tau_0$. Notice that when $k<\frac{n-2}{2}$, $k<\tau_0<\frac{n-2}{2}$. Therefore the set for $\tau$ is not empty when $n\ge 6$. When $n=5$, $k=1$ since $k<\frac{n-2}{2}$. In this case $n-2<\frac{n-2}{2}+2+\tau-\frac{4}{n-2}k$ and we can just choose $k\le \tau<\tau_0=\frac{n-2}{2}$.

\medskip

When $y\in \Omega_i\cap B_i$ for some $i$,
$$
\begin{array}{l}
\int_{\cup_j(\Omega_j\cap B_j^c)}\frac{1}{|y-z|^{n-2}}W_m^{\frac{4}{n-2}}\gamma(z)\sum_s\frac{1}{(1+|z-X^s|)^{\frac{n-2}{2}+\tau}}dz \\
\\
\le C\frac{1}{(\lambda l)^{\frac{4}{n-2}k}}\int_{{\mathbb R}^n}\frac{1}{|y-z|^{n-2}}\frac{1}{\lambda^{\tau-1}}\sum_j\frac{1}{(1+|z-X^s|)^{4+\frac{n-2}{2}+1-\frac{4}{n-2}k}}dz\\
\\
\le C\frac{1}{\lambda^{\tau-1}}\sum_j\frac{1}{(1+|y-X^j|)^{\min(\frac{n-2}{2}+2+1-\frac{4}{n-2}k,n-2)}}\frac{1}{(\lambda l)^{\frac{4}{n-2}k}}\\
\\
\le \frac{C}{\lambda^{\tau-1}}\frac{1}{(1+|y-X^i|)^{\min(\frac{n-2}{2}+2+1-\frac{4}{n-2}k,n-2)}}\frac{1}{(\lambda l)^{\frac{4}{n-2}k}},\quad\mbox{due to}\quad y\in \Omega_i\cap B_i\\
\\
\le \frac{C}{(\lambda l)^{\frac{4k}{n-2}}}\left\{\begin{array}{l}
\frac{\gamma(y)}{(1+|y-X^i|)^{\frac{n-2}{2}+\tau}},\quad\mbox{if $n-2> \frac{n-2}{2}+3-\frac{4}{n-2}$}\\
\\
\frac{\gamma(y)}{(1+|y-X^i|)^{n-2+\tau-1}}\le \frac{\gamma(y)}{(1+|y-X^i|)^{\frac{n-2}{2}+\tau}}\quad\mbox{if $n-2< \frac{n-2}{2}+3-\frac{4}{n-2}$}
\\
\end{array}\right.\\
\\
\le C\gamma(y)\sum_j\frac{1}{(1+|y-X^j|)^{\frac{n-2}{2}+\tau}}\frac{1}{(\lambda l)^{\frac{4}{n-2}k}}.
\end{array}
$$

\medskip

When $\frac{n-2}{2}+2+1-\frac{4}{n-2}k = n-2$, the $\log |y|$ term from applying Lemma \ref{L2} will be harmless as long as  we choose $\tau<\tau_0\le \frac{n-2}{2}$. The fact that $2>\frac{4k}{n-2}$ is also used in the above. Combining the  above together, for $0<\theta<\min(\theta_1,\theta_2,\theta_3)$, we conclude that
$$\begin{array}{c}
\int_{{\mathbb R}^n}\frac{1}{|y-z|^{n-2}}W_m^{\frac{4}{n-2}}\gamma(z)\sum_j\frac{1}{(1+|z-X^j|)^{\frac{n-2}{2}+\tau}}dz\le\\
\\
\sum_j\frac{C\gamma(y)}{(1+|y-X^j|)^{\frac{n-2}{2}+\theta+\tau}}+\frac{C\gamma(y)}{(\lambda l)^{\frac{4}{n-2}k}}\sum_j\frac{1}{(1+|y-X^j|)^{\frac{n-2}{2}+\tau}}.
\end{array}
$$
\end{proof}

\medskip

\renewcommand{\theLemma}{A.5}
\begin{Lemma}\label{L8}
Assume  $n\ge 4$ and $0< \tau<\frac{n+2}{2}$. If  $\|\phi\|_*\le \frac{C}{\lambda^{\frac{n+2}{2}-\tau}}$,  then for any $c>0$, there exists a constant $\lambda_0=\lambda_0(n,k,\tau,C,c)>0$, such that for any $\lambda>\lambda_0$, $\phi(y)\le cW_m(y)$ in $\cup_i (B_i\cap \Omega_i)$.
\end{Lemma}

\begin{proof}
We prove by contradiction. Without loss of generality, we may assume that $\phi(y)\ge cW_m(y)$ for some $y\in B_1\cap\Omega_1$. Note that $\gamma(y)\le 1$. By Lemma \ref{L5}, we have
$$\begin{array}{c}
\phi(y)\ge cC\sum_j\frac{1}{(1+|y-X^j|)^{n-2}}\ge C\frac{1}{(1+|y-X^1|)^{n-2}}\\
\\
\ge C\frac{1}{(1+|y-X^1|)^{\frac{n-2}{2}+\tau}}\frac{1}{(1+|y-X^1|)^{\frac{n-2}{2}-\tau}}
\\
\\
\ge  C\gamma(y)\sum_j\frac{1}{(1+|y-X^j|)^{\frac{n-2}{2}+\tau}}\frac{1}{(1+|y-X^1|)^{\frac{n-2}{2}-\tau}}.
\end{array}
$$
When $n\ge 4$ and $0< \tau<\frac{n-2}{2}$, we have
$$
\frac{C}{\lambda^{\frac{n+2}{2}-\tau}}\ge\|\phi\|_*\ge \frac{1}{(\lambda l)^{\frac{n-2}{2}-\tau}}.
$$

This gives a contradiction when $\lambda$ is large. In the case $n\ge 4$ and $\frac{n+2}{2}>\tau\ge \frac{n-2}{2}$, noting the fact that
$$
\frac{1}{(1+|y-X^1|)^{\frac{n-2}{2}-\tau}}\ge 1, \quad \forall y\in B_1,
$$we get $\frac{C}{\lambda^{\frac{n+2}{2}-\tau}}\ge\|\phi\|_*\ge C$, which is impossible when $\lambda$ large.
\end{proof}

\medskip

\renewcommand{\theLemma}{A.6}
\begin{Lemma}\label{La.12}
For any $\phi\in\mathcal{\tilde M}$, we have, for some $C>0$, independent of $m$ and $l$,
$$
|\int_{{\mathbb R}^n}K_{\lambda}(z)W_m^{\frac{4}{n-2}}\phi Z_{s,t}dz|\le \frac{C\|\phi\|_*}{\lambda^{\frac{n-2}{2}+\tau}},
$$
and
$$
|\int_{{\mathbb R}^n}K_{\lambda}(z)W_m^{\frac{4}{n-2}}\phi \sigma_idz|\le \frac{C\|\phi\|_*}{\lambda^{\frac{n-2}{2}+\tau}}.
$$
\end{Lemma}
\begin{proof} By the orthogonality condition
$$
\begin{array}{l}
\int_{{\mathbb R}^n}K_{\lambda}(z)W_m^{\frac{4}{n-2}}\phi Z_{s,t}dz=\int_{{\mathbb R}^n}(K_{\lambda}(z)-1)\sigma_s^{\frac{4}{n-2}}Z_{s,t}\phi dz\\
\\
+O(1)|\int_{{\mathbb R}^n}\hat{W}_{m,s}^{\frac{4}{n-2}}Z_{s,t}\phi dz|+O(1)|\int_{{\mathbb R}^n}\hat{W}_m\sigma_s^{\frac{4}{n-2}-1}Z_{s,t}\phi dz|.
\end{array}
$$

\medskip

Using Lemma \ref{L5} and the proof of Lemma \ref{L3} in $\Omega_i$ with $i\neq s$, we get
$$
\begin{array}{l}
|\int_{\Omega_i\cap B_i}\hat{W}_{m,s}^{\frac{4}{n-2}}Z_{s,t}\phi dz|\le C\|\phi\|_*\int_{\Omega_i\cap B_i}\sum_j\frac{\gamma(z)}{(1+|z-X^j|)^{\frac{n-2}{2}+\tau}}\frac{1}{(1+|z-X^s|)^{n-2}}\times\\
\\
\left(\sum_{k\neq s}\frac{1}{(1+|z-X^k|)^{n-2}}\right)^{\frac{4}{n-2}}dz\\
\\
\le \frac{C\|\phi\|_*}{\lambda^{\tau-1}}\int_{\Omega_i\cap B_i}\frac{1}{(1+|z-X^i|)^{\frac{n-2}{2}+5}}\frac{1}{(1+|z-X^s|)^{n-2}}dz\\
\\
\le \frac{C\|\phi\|_*}{\lambda^{\tau-1}|X^i-X^s|^{\frac{n}{2}}}.
\end{array}
$$

When $z\in \cup_i(\Omega_i\cap B_i^c)$, we use the same idea as in the proof for Lemma \ref{L3} to get
$$
\begin{array}{l}
|\int_{\cup_i(\Omega_i\cap B_i^c)}\hat{W}_{m,s}^{\frac{4}{n-2}}Z_{s,t}\phi dz|\\
\\
\le \frac{C\|\phi\|_*}{(\lambda l)^{\frac{4}{n-2}k}}\int_{\cup_i(\Omega_i\cap B_i^c)}\sum_j\frac{1}{(1+|z-X^j|)^{\frac{n-2}{2}+\tau+4-\frac{4}{n-2}k}}\frac{1}{(1+|z-X^s|)^{n-2}}dz\\

\\
\le \frac{C\|\phi\|_*}{(\lambda l)^{\frac{4}{n-2}k}} \left(\sum_{j\neq s}\frac{1}{(1+|X^j-X^s|)^{\min(n-2,\frac{n-2}{2}+\tau+2-\frac{4}{n-2}k)}}+\frac{1}{(\lambda l)^{2+\frac{n-2}{2}+\tau-\frac{4}{n-2}k}}\right)\\
\\
\le \frac{C\|\phi\|_*}{\lambda^{\tau-1}(\lambda l)^{\frac{n}{2}}}.
\end{array}
$$

For $i=s$, note that for any $z\in \Omega_s$, taking $X^j$ be the closest point in $X_{l,m}$ to $X^s$ (there are at most $2^k$ such kind of points in $X_{l,m}$), we have
$$
\hat{W}_{m,s}^{\frac{4}{n-2}}\le\frac{ C}{(1+|z-X^j|)^{4-\frac{4}{n-2}k}(\lambda l)^{\frac{4}{n-2}k}}.$$

By Lemma \ref{L1} we also get
$$
|\int_{\Omega_s\cap B_s}\hat{W}_{m,s}^{\frac{4}{n-2}}Z_{s,t}\phi dz|\le  \frac{C\|\phi\|_*}{(\lambda l)^{\frac{n}{2}}\lambda^{\tau-1}}.
$$

Hence we deduce that
$$
|\int_{{\mathbb R}^n}\hat{W}_{m,s}^{\frac{4}{n-2}}Z_{s,t}\phi dz|\le \frac{C\|\phi\|_*}{(\lambda l)^{\frac{n}{2}} \lambda^{\tau-1}}
$$for all $n\ge 5$ and $k\le\tau<\tau_0$.

\medskip

Similarly we have
$$\begin{array}{l}
|\int_{{\mathbb R}^n}\hat{W}_{m,s}\sigma_s^{\frac{4}{n-2}-1}Z_{s,t}\phi dz|\le C\int_{{\mathbb R}^n}\hat{W}_{m,s}\sigma_s^{\frac{4}{n-2}}|\phi|dz\\
\\
\le \frac{C\|\phi\|_*}{(\lambda l)^{\frac{n}{2}}\lambda^{\tau-1}}.
\end{array}
$$

$$\begin{array}{l}
|\int_{{\mathbb R}^n}(K_{\lambda}(z)-1)\sigma_s^{\frac{4}{n-2}}Z_{s,t}\phi dz|\\
\\
\le C\|\phi\|_*\int_{{\mathbb R}^n}|K_{\lambda}(z)-1| \frac{\gamma(z)}{(1+|z-X^s|)^{n+2}}\sum_j \frac{1}{(1+|z-X^j|)^{\frac{n-2}{2}+\tau}}dz.
\end{array}
$$

It is easy to see by Lemma \ref{L1}  and the relations between $X^i$ that
$$
|\int_{{\mathbb R}^n}\frac{|K_{\lambda}(z)-1| \gamma(z)}{(1+|z-X^s|)^{n+2}}\sum_{j\neq s} \frac{1}{(1+|z-X^j|)^{\frac{n-2}{2}+\tau}}dz|\le\frac{C}{(\lambda l)^{\frac{n-2}{2}+\tau}}.
$$

So we only need to estimate the integral when $j=s$. Let $\Omega=\{z\in{\mathbb R}^n| |z-X^s|\le \lambda$. By Lemma \ref{L5} and integrating, it is not hard to get
$$\begin{array}{l}
\int_{\Omega}|K_{\lambda}(z)-1| \frac{\gamma(z)}{(1+|z-X^s|)^{n+2}} \frac{1}{(1+|z-X^s|)^{\frac{n-2}{2}+\tau}}dz\\
\\
\le\frac{C}{\lambda^{\beta+\tau-1}}\int_{\Omega}\frac{|z-X^s|^{\beta+\tau-1}}{(1+|z-X^s|)^{n+2+\frac{n-2}{2}+\tau}}dz\\
\\
\le C\left\{\begin{array}{l}
\frac{1}{\lambda^{\frac{n+2}{2}+\tau}}\quad\mbox{if}\quad \beta>\frac{n}{2}+2,\\
\\
\frac{1}{\lambda^{\beta+\tau-1}}\le \frac{1}{\lambda^{\frac{n-2}{2}+\tau}}\quad\mbox{if} \quad \beta\le \frac{n}{2}+2,
\end{array}\right.
\end{array}
$$
$$\begin{array}{l}
\int_{\Omega^c}|K_{\lambda}(z)-1| \frac{\gamma(z)}{(1+|z-X^s|)^{n+2}}\frac{1}{(1+|z-X^s|)^{\frac{n-2}{2}+\tau}}dz\\
\\
\le C\int_{\Omega^c}\frac{1}{(1+|z-X^s|)^{n+2+\frac{n-2}{2}+\tau}}dz\\
\\
\le \frac{C}{\lambda^{\frac{n+2}{2}+\tau}}.
\end{array}
$$
 since $\beta>n-2>\frac{n}{2}$ when $n\ge 5$.

Therefore the first inequality can be derived easily and the second one can be proved similarly.
\end{proof}

\medskip

\renewcommand{\theLemma}{A.7}
\begin{Lemma}\label{La.11}
Under the assumption of Lemma \ref{L3}, for any $h\in\mathcal{\tilde D}$ and $\phi\in\mathcal{\tilde{M}}$, let $\tilde{\phi}=P_{\mathcal E}(-\Delta)^{-1}(h+K_{\lambda}W_m^{\frac{4}{n-2}}\phi)$, then there exist an integer $l_0$ and a constant $C>0$, depending only on $K$, $n$, $\beta$, $\tau$, $C_1$ and $C_2$, such that for any $l\ge l_0$, we have
$$
\|\tilde{\phi}\|_*\le C(\|h\|_{**}+\|\phi\|_*).
$$
\end{Lemma}

\begin{proof}
By assumption on $\phi$, $\phi$ satisfies the equation
\begin{equation*}
\tilde{\phi}(y)=\frac{1}{n(n-2)\omega_n}\int_{{\mathbb R}^n}\frac{h+K_{\lambda}W_m^{\frac{4}{n-2}}\phi}{|y-z|^{n-2}}dz+\sum_{i,j}c_{i,j}Z_{i,j}+\sum_i b_i\sigma_i,
\end{equation*}for some constants $c_{i,j}$, $b_i$.

\medskip

We first claim that, for some constant $C$, independent of $m$ and $l$,
\begin{equation*}\tag{A6}\label{A7.1}
|c_{i,j}| , |b_i| \le\left( C\|h\|_{**}+\frac{C}{\lambda^{\frac{n}{2}}}\|\phi\|_*\right)\frac{1}{\lambda^{\tau-1}}.
\end{equation*}

In fact, multiplying $\sigma_s^{\frac{4}{n-2}}Z_{s,t}$ on both side of the equation and integrating, we get
\begin{equation*}\tag{A7}\label{A8.1}
\sum_{i,j}c_{i,j}\l Z_{i,j},Z_{s,t}\r=\int_{{\mathbb R}^n}\left(-h-K_{\lambda}W_m^{\frac{4}{n-2}}\phi-\sum_j b_i\sigma_i^{\frac{n+2}{n-2}}\right)Z_{s,t}dz,
\end{equation*}

$$
\begin{array}{c}
|\int_{{\mathbb R}^n}h(z)Z_{s,t}dz|\le C\|h\|_{**}\int_{{\mathbb R}^n}\frac{1}{(1+|z-X^s|)^{n-2}}\gamma(z)\sum_j\frac{1}{(1+|z-X^j|)^{\frac{n+2}{2}+\tau}}dz\\
\\
\le C\|h\|_{**}\left(\int_{{\mathbb R}^n}\frac{\gamma(z)}{(1+|z-X^s|)^{n-2+\frac{n+2}{2}+\tau}}dz\right.\\
\\
\left. +\sum_{j\neq s}\int_{{\mathbb R}^n}\frac{\gamma(z)}{(1+|z-X^s|)^{n-2}}\frac{1}{(1+|z-X^j|)^{\frac{n+2}{2}+\tau}}dz\right),
\end{array}
$$
where
$$\begin{array}{l}
\int_{{\mathbb R}^n}\frac{\gamma(z)}{(1+|z-X^s|)^{n-2+\frac{n+2}{2}+\tau}}dz\le\int_{B_s}\frac{1}{\lambda^{\tau-1}}
\frac{1}{(1+|z-X^s|)^{n-2+\frac{n+2}{2}+1}}dz\\
\\
+\int_{B_s^c}\frac{1}{(1+|z-X^s|)^{n-2+\frac{n+2}{2}+\tau}}dz\\
\\
\le \frac{C}{\lambda^{\tau-1}}+\frac{C}{(\lambda l)^{\frac{n-2}{2}+\tau}}\le \frac{C}{\lambda^{\tau-1}},
\end{array}
$$
and
$$\begin{array}{l}
\sum_{j\neq s}\int_{{\mathbb R}^n}\frac{\gamma(z)}{(1+|z-X^s|)^{n-2}}\frac{1}{(1+|z-X^j|)^{\frac{n+2}{2}+\tau}}dz\\
\\
\le \frac{1}{\lambda^{\tau-1}}\sum_{j\neq s}\int_{{\mathbb R}^n}\frac{1}{(1+|z-X^s|)^{n-2-\tau+1}}\frac{1}{(1+|z-X^j|)^{\frac{n+2}{2}+\tau}}dz\\
\\
\le \frac{1}{\lambda^{\tau-1}}\sum_{j\neq s}\frac{1}{|X^s-X^j|^{\frac{n-2}{2}}}(\int_{{\mathbb R}^n}\frac{1}{(1+|z-X^j|)^{n+1}}+\frac{1}{(1+|z-X^s|)^{n+1}}dz)\\
\\
\le \frac{C}{\lambda^{\tau-1}}.
\end{array}
$$

Here we have used the fact that $\sum_{j\neq s}\frac{1}{|X^j-X^s|^{\frac{n-2}{2}}}$ converges when $1\le k<\frac{n-2}{2}$. Thus we have derived
$$
|\int_{{\mathbb R}^n}h(z)Z_{s,t}dz|\le \frac{C}{\lambda^{\tau-1}}\|h\|_{**}.
$$

By Lemma \ref{La.12},
$$
|\int_{{\mathbb R}^n}K_{\lambda}(z)W_m^{\frac{4}{n-2}}\phi Z_{s,t}dz|\le \frac{C\|\phi\|_*}{\lambda^{\frac{n-2}{2}+\tau}},
$$
By Lemma \ref{L1} and symmetry of $\sigma_i$, it is easy to check that,
$$\l Z_{i,j},Z_{s,t}\r=0 \quad\mbox{if} \quad i=s\quad\mbox{and}\quad j\neq t;$$

$$
\l Z_{i,j},Z_{i,j}\r=C;
$$
and

$$|\l Z_{i,j},Z_{s,t}\r|\le \frac{C}{|X^i-X^s|^{n-2}},\quad\mbox{if}\quad i\neq s,$$

$$\int_{{\mathbb R}^n}\sigma_i^{\frac{n+2}{n-2}}Z_{s,t}=0,\quad\mbox{if}\quad i=s;$$
and
$$
| \int_{{\mathbb R}^n}\sigma_i^{\frac{n+2}{n-2}}Z_{s,t} | \le \frac{C}{|X^i-X^s|^{n-2}}, \quad\mbox{if}\quad i\neq s.
$$

Notice that the left hand of equation (\ref{A8.1}) can be viewed as a linear system with variables $c_{i,j}$ of $(m+1)^k(n+1)$ dimension and coefficient matrix of $(m+1)^k(n+1)\times (m+1)^k(n+1)$ with entry $\l Z_{i,j},Z_{s,t}\r$. If we denote $G=\left(\l Z_{i,j}, Z_{s,t}\r\right)=\left(a_{i,j}^{s,t}\right)$  be this matrix and let $X=(x_{s,t})$ be in ${\mathbb R}^{(m+1)^k(n+1)}$ with maximum norm denotes as $|X|=\max_{s,t}|x_{s,t}|$, then
$$\begin{array}{l}
C|x_{i,j}|+\frac{c(n+1)}{(\lambda l)^{n-2}}|X|\ge |\sum_{s,t}a_{i,j}^{s,t}x_{s,t}|\\
\\
=|Cx_{i,j}+\sum_{(s,t)\neq(i,j)}a_{i,j}^{s,t}x_{s,t}|\ge C|x_{i,j}|-\frac{c(n+1)}{(\lambda l)^{n-2}}|X|,
\end{array}
$$
where $c$ is controlled by $\int_{{\mathbb R}^k}\frac{1}{1+|z|^{n-2}}dz$ and doesn't depend on $m$. This implies that
$$
2C|X|\ge |GX|\ge \frac{C}{2}|X|
$$
with $C$ independent of $m$ when $\lambda$ is large enough. Therefore, we obtain that
\begin{equation*}
|c_{i,j}|\le \left( C\|h\|_{**}+\frac{C}{\lambda^{\frac{n}{2}}}\|\phi\|_*\right)\frac{1}{\lambda^{\tau-1}}+C\max|b_i|\frac{1}{(\lambda l)^{n-2}}.
\end{equation*}where $C$ is a constant that doesn't depend on $m$ and $l$.

Similar estimates for $b_i$ can be obtained in the same way as $c_{i,j}$, we skip the detail. Thus we prove (\ref{A7.1}).

\medskip

It follows from Lemma \ref{L2},
$$
\begin{array}{c}
|\int_{{\mathbb R}^n}\frac{1}{|y-z|^{n-2}}h(z)dz|\le C\|h\|_{**}\int_{{\mathbb R}^n}\frac{\gamma(z)}{|y-z|^{n-2}}\sum_j\frac{1}{(1+|z-X^j|)^{\frac{n+2}{2}+\tau}}dz\\
\\
\le C\|h\|_{**}\gamma(y)\sum_j\frac{1}{(1+|y-X^j|)^{\frac{n-2}{2}+\tau}},
\end{array}
$$

and
$$
\begin{array}{c}
|\int_{{\mathbb R}^n}\frac{1}{|y-z|^{n-2}}\sigma_i^{\frac{4}{n-2}}Z_{i,j}dz|\\
\\
\le C \int_{{\mathbb R}^n}\frac{1}{|y-z|^{n-2}}\frac{1}{(1+|z-X^i|)^{n+2}}dz\le \frac{C}{(1+|y-X^i|)^{\frac{n-2}{2}+\tau}},\quad \forall k\le\tau< \tau_0.
\end{array}
$$

Similarly,
$$
|\int_{{\mathbb R}^n}\frac{1}{|y-z|^{n-2}}\sigma_i^{\frac{n+2}{n-2}}dz|
\le \frac{C}{(1+|y-X^i|)^{\frac{n-2}{2}+\tau}}
$$
when $k\le \tau< \tau_0$.

This, combined with Lemma \ref{L3} and (\ref{A7.1}), gives the conclusion. The proof of Lemma is thus completed.
\end{proof}

\medskip

\renewcommand{\theLemma}{A.8}
\begin{Lemma}\label{La.13}
Under the same assumption of Lemma \ref{L3}, there exist an integer $l_0$ and a constant $C\ge 1$, depending only on $K$, $n$, $\beta$, $\tau$, $C_1$ and $C_2$, such that for any $\phi\in\mathcal{\tilde E}$, we have
\begin{equation*}\tag{A8}\label{A9.1}
\|\phi-\frac{n+2}{n-2}P_{\mathcal E}(-\Delta)^{-1}\left(K_{\lambda}W_m^{\frac{4}{n-2}}\phi\right)\|_*\ge \frac{\|\phi\|_*}{C}.
\end{equation*}
\end{Lemma}
\begin{proof}
(\ref{A9.1}) is equivalent to
\begin{equation*}\tag{A9}\label{A10.1}
\begin{array}{l}
\phi(y)=h+\frac{n+2}{n(n-2)^2\omega_n}\int_{{\mathbb R}^n}\frac{1}{|z-y|^{n-2}}K_{\lambda}(z)W_m^{\frac{4}{n-2}}(z)\phi(z)dz\\
\\
+\sum_i b_i\sigma_i+\sum_{i,j}c_{i,j}Z_{i,j},
\end{array}
\end{equation*}for some $h\in{\mathcal{\tilde E}}$ and constants $b_i$, $c_{i,j}$.  Then by Lemma \ref{L3} and the proof of Lemma \ref{La.11}, we get
\begin{equation*}\tag{A10}\label{A11.1}
\begin{array}{l}
\left(\gamma(y)\sum_{X^i\in X_{l,m}}\frac{1}{(1+|y-X^i|)^{\frac{n-2}{2}+\tau}}\right)^{-1}|\phi(y)\le \|h\|_*+\frac{C\|\phi\|_*}{\lambda^{\frac{n-2}{2}+\tau}}\\
\\
+C\|\phi\|_*\left(\frac{1}{(\lambda l)^{\frac{4k}{n-2}}}+\frac{\sum_{X^j\in X_{l,m}}\frac{1}{(1+|y-X^j|)^{\frac{n-2}{2}+\tau+\theta}}}{\sum_{X^j\in X_{l,m}}\frac{1}{(1+|y-X^j|)^{\frac{n-2}{2}+\tau}}}\right).
\end{array}
\end{equation*}

\medskip

We show that $\|\phi\|_*\le C\|h\|_*$ for $l$ large enough. If not, we can find sequences $l\to\infty$, $m_l\ge 1$, $\Lambda_{i,l}\in [C_1, C_2]$, $P^{i,l}\in B_{\frac{1}{2}}(X^i)$, $b_{i,l}$, $c_{i,j,l}$ and $\phi_l$ with $\|\phi_l\|_*=1$, such that (\ref{A10.1}) is satisfied for $\|h_l\|_*\to 0$. We may assume that $\|\phi_l\|_*=1$. Therefore for some $y^l\in {\mathbb R}^n$, we obtain from (\ref{A11.1}) that
 \begin{equation*}\label{A12.1}\tag{A11}
1=\|\phi_l\|_*\le C\left(\|h_l\|_{*}+\frac{\sum_{X^j\in X_{l,m}}\frac{1}{(1+|y^l-X^j|)^{\frac{n-2}{2}+\tau+\theta}}}{\sum_{X^j\in X_{l,m}}\frac{1}{(1+|y^l-X^j|)^{\frac{n-2}{2}+\tau}}}\right).
\end{equation*}

\medskip

Then there is a $R>0$ independent of $m$, such that
for some $i(l)$, $y^l\in B_R(X^{i(l)})$ for all $l$ large. (If $y^l$ is far away from all $X^i$, the right side of (\ref{A12.1}) is approaching zero as $l\to\infty$.) Hence we get that
 $$
\max_{{\mathbb B}_R(x^{i(l)})} |\lambda^{\tau-1}\phi_l(y)|\ge a>0.
$$

From the proof of Lemma \ref{L5}, for any fixed $R>0$, it is easy to see that $W_m(x-P^{i(l)})\to \sigma_{0,\Lambda}$ for some $\Lambda\in [C_1,C_2]$ in $B_R$ as $l\to\infty$ independent of $m$. Multiplying $\lambda^{\tau-1}$ on both side of the equation (\ref{A10.1}) and using the estimates (\ref{A7.1}) for $b_{i,l}$ and $c_{i,j,l}$ and  the fact that $\|h_l\|_*\to 0$, we can see that $\tilde{\phi}_l(y):=\lambda^{\tau-1}\phi(y-P^{i(l)})$ converges uniformly in any compact set to a non-zero solution $\tilde{\phi}$ of
\begin{equation*}\tag{A12}\label{A13.1}
-\Delta \tilde{\phi}-\frac{n+2}{n-2}\sigma_{0,\Lambda}^{\frac{4}{n-2}}\tilde{\phi}=0
\end{equation*}
for some $\Lambda\in [C_1,C_2]$. We will show that $\tilde{\phi}$ is perpendicular to the kernel of (\ref{A13.1}) and therefore $\tilde{\phi}=0$, which is a contradiction.

\medskip

For this purpose, by Lemma \ref{L5}, we get
$$
\lambda^{\tau-1}\phi(y)\le \Phi(y)=C\left\{\begin{array}{l}
                          \frac{1}{(1+|y-X^j|)^{\frac{n-2}{2}+1}}\quad y\in B_j\cap\Omega_j\\
\\
\frac{\lambda^{\tau-1}}{(\lambda l)^k}\frac{1}{(1+|y-X^j|)^{\frac{n-2}{2}+\tau-k}}\quad y\in B_j^c\cap\Omega_j.
                          \end{array}\right.
$$

Since
$$
\sigma_i(y),\quad |\frac{\p \sigma_i}{\p {P^i}_j}(y)|,\quad |\frac{\p \sigma_i}{\p \Lambda_i}(y)|\le \frac{C}{(1+|y-X^i|)^{n-2}},
$$
from dominated convergence theorem and the fact that $\l\Phi,\sigma_i\r=0,\quad \l\Phi,\frac{\p\sigma_i}{\p {P^i}_j}\r=0$ and $\l\Phi,\frac{\p\sigma_i}{\p \Lambda_i}\r=0$, we obtain
$$
\l\tilde{\phi},\sigma_{0,\Lambda}\r,\quad \l\tilde{\phi},\frac{\p\sigma_{0,\Lambda}}{\p x_j}\r,\quad\l\tilde{\phi},\frac{\p\sigma_{0,\Lambda}}{\p\Lambda}\r=0.
$$
Therefore we have proved the conclusion.
\end{proof}

\medskip

\renewcommand{\theLemma}{A.9}
\begin{Lemma}\label{La.1}
For $j=1,...,n$, $i=1,...,(m+1)^k$, $0<C_1\le\Lambda_i\le C_2<\infty$ and $P^i\in B_{\frac{1}{2}}(X^i)$, we have
$$
\begin{array}{l}
\int_{{\mathbb R}^n}
K_{\lambda}(y)\sigma_i^{\frac{n+2}{n-2}}\frac{\partial\sigma_i}{\partial
{P^i}_j}=\frac{D_{n,\beta}a_j}{\Lambda_i^{\beta-2}\lambda^{\beta}}({P^i}_j-{X^i}_j)\\
\\
+O(\frac{|P^i-X^i|^2}{\lambda^{\beta}})+o(\frac{1}{\lambda^{\beta}}).
\end{array}
$$where $D_{n,\beta}=(n(n-2))^{\frac{n}{2}}(n-2)\beta\int_{{\mathbb R}^n}\frac{|x_j|^{\beta}}{(1+|x|^2)^{n+1}}dx$.  $o_{\lambda}(1)$ only depends on the condition of function $R(\frac{y}{\lambda})$ near $X^i$ and $o(1)\to 0$ as $l\to\infty$ (or same as $\lambda\to \infty$) (see the remark \ref{R2} below).
\end{Lemma}

\medskip

\begin{proof}
Let $\delta=\lambda^{\frac{\beta-n}{2n}}$. We have
$$
\begin{array}{l}
\int_{{\mathbb R}^n}
K_{\lambda}(y)\sigma_i^{\frac{n+2}{n-2}}\frac{\partial\sigma_i}{\partial
{P^i}_j}=(n-2)\int_{{\mathbb R}^n}K_{\lambda}(y)\sigma_i^{\frac{2n}{n-2}}\frac{\Lambda_i^2(y_j-{P^i}_j)}{1+\Lambda_i^2|y-P^i|^2}dy\\
\\
=(n-2)\int_{|y-X^i|\le\delta\lambda}K_{\lambda}(y)\sigma_i^{\frac{2n}{n-2}}\frac{\Lambda_i^2(y_j-{P^i}_j)}{1+\Lambda_i^2|y-P^i|^2}dy+O(\frac{1}{(\delta\lambda)^{n+1}})\\
\\
=\frac{n-2}{\lambda^{\beta}}\int_{|y-X^i|\le\delta\lambda}\left(\sum_h a_h|y_h-{X^i}_h|^{\beta}+o(1)|y-X^i|^{\beta}\right)\\
\\
\times \sigma_i^{\frac{2n}{n-2}}\frac{\Lambda_i^2(y_j-{P^i}_j)}{1+\Lambda_i^2|y-P^i|^2}dy
+O(\frac{1}{(\delta\lambda)^{n+1}})\\
\\
=\frac{n-2}{\lambda^{\beta}}\int_{|y-X^i|\le\delta\lambda}\sum_h a_h|y_h-{X^i}_h|^{\beta}\sigma_i^{\frac{2n}{n-2}}\frac{\Lambda_i^2(y_j-{P^i}_j)}{1+\Lambda_i^2|y-P^i|^2}dy\\
\\
+o(\frac{1}{\lambda^{\beta}})+o(\frac{|P^i-X^i|^{\beta}}{\lambda^{\beta}})\\
\\
=(n(n-2))^{\frac{n}{2}}\frac{n-2}{\lambda^{\beta}}\int_{{\mathbb R}^n}\sum_ha_h(|x_h|^{\beta}+\beta|x_h|^{\beta-2}x_h(P^i-X^i)_h\\
\\
+O(|P^i-X^i|^2))\times \frac{\Lambda_i^n}{(1+\Lambda_i^2|x|^2)^n}\frac{\Lambda_i^2 x_j}{(1+\Lambda_i^2|x|^2)}dx+o(\frac{1}{\lambda^{\beta}})+o(\frac{|P^i-X^i|^{\beta}}{\lambda^{\beta}})\\
\\
=(n(n-2))^{\frac{n}{2}}\frac{(n-2)\beta a_j}{\Lambda_i^{\beta-2}\lambda^{\beta}}\int_{{\mathbb R}^n}\frac{|x_j|^{\beta}}{(1+|x|^2)^{n+1}} ({P^i}_j-{X^i}_j)\\
\\
+o(\frac{1}{\lambda^{\beta}})+O(\frac{|P^i-X^i|^2}{\lambda^{\beta}})+o(\frac{|P^i-X^i|^{\beta}}{\lambda^{\beta}}).
\end{array}
$$

If we let $D_{n,\beta}=(n(n-2))^{\frac{n}{2}}(n-2)\int_{{\mathbb R}^n}\frac{|x_j|^{\beta}}{(1+|x|^2)^{n+1}}dx$, we complete the proof.
\end{proof}

\medskip

\begin{remark} \em \label{R2}
 The $o_\lambda(1)$ only depends on the condition of $R(\frac{y}{\lambda})$ near $X^i$, the estimates doesn't depend on $P_i$ as long as $|P_i-X^i|\le \frac{1}{2}$.  If we know more, say $|\nabla R(x)\le C|x|^{\beta-1+s}$ near $0$ for some small $s>0$, then $o_\lambda(1)=\frac{C}{\lambda^s}$.
\end{remark}

\medskip

\renewcommand{\theLemma}{A.10}
\begin{Lemma}\label{La.3}
 $$
\begin{array}{l}
\int_{{\mathbb R}^n}
K_{\lambda}(y)\sigma_i^{\frac{n+2}{n-2}}\frac{\partial\sigma_i}{\partial
\Lambda_i}=\frac{C_3}{\Lambda_i^{\beta+1}\lambda^{\beta}}\\
\\
+O(\frac{|P^i-X^i|^{\beta-1}}{\lambda^{\beta}})+o(\frac{1}{\lambda^{\beta}}).
\end{array}
$$where $o(1)$ is the same as in Lemma \ref{La.1} and
$$
C_3=-\frac{\beta[n(n-2)]^{\frac{n}{2}}(n-2)}{2n}(\sum_i a_i)\int_{{\mathbb R}^n}\frac{|y_1|^{\beta}}{(1+|y|^2)}dy>0.
$$
\end{Lemma}
\begin{proof} Observe that
$$
\begin{array}{l}
\int
K_{\lambda}(y)\sigma_i^{\frac{n+2}{n-2}}\frac{\partial\sigma_i}{\partial
\Lambda_i}=\frac{n-2}{2n}\frac{\partial}{\partial \Lambda_i}\int
K_{\lambda}(y)\sigma_i^{\frac{2n}{n-2}}\\
\\
=-\frac{1}{\Lambda_i}\frac{[n(n-2)]^{\frac{n}{2}}(n-2)}{2n}\int_{|y|\le\delta\lambda}(\nabla
K(\frac{y}{\lambda\Lambda_i}+\frac{P^i}{\lambda})\cdot \frac{y}{\lambda\Lambda_i})\frac{1}{(1+|y|^2)^n}dy+O(\frac{1}{(\delta \lambda)^n})\\
\\
=-\frac{\beta}{\Lambda_i}\frac{[n(n-2)]^{\frac{n}{2}}(n-2)}{2n}\int_{{\mathbb R}^n}\sum a_h\frac{|y_h|^{\beta}}{(\Lambda_i\lambda)^{\beta}}\frac{1}{(1+|y|^2)^n}dy+\\
\\
+o(\frac{1}{\lambda^{\beta}})+O(\frac{|P^i-X^i|^{\beta-1}}{\lambda^{\beta}})\\
\\
=\frac{C_3}{\Lambda_i^{\beta+1}\lambda^{\beta}}+o(\frac{1}{\lambda^{\beta}})+O(\frac{|P^i-X^i|^{\beta-1}}{\lambda^{\beta}}).
\end{array}
$$
\end{proof}

\medskip

\renewcommand{\theLemma}{A.11}
\begin{Lemma}\label{La.2}
For $j=1,...,n$, we have
$$\begin{array}{l}
\int_{{\mathbb R}^n}K_{\lambda}(y)(\bar{W}_m+\phi)^{\frac{n+2}{n-2}}Z_{i,j}dy=\int_{{\mathbb R}^n}K_{\lambda}(y)\sigma_i^{\frac{n+2}{n-2}}Z_{i,j}dy\\
\\
+C\left(|\epsilon|^2+\|\phi\|_*^{\frac{n+2}{n-2}}\frac{1}{(\lambda l)^{\frac{n}{2}+\tau\frac{n+2}{n-2}}}+\|\phi\|_*^2\frac{1}{\lambda^{2\tau-2}}+\frac{1}{\lambda^{\beta+1}}+o(\frac{1}{\lambda^n})\right).
\end{array}
$$
\end{Lemma}

\begin{proof} We begin with
$$
\begin{array}{c}
|(\bar{W}_m+\phi)^{\frac{n+2}{n-2}}-W_m^{\frac{n+2}{n-2}}-\frac{n+2}{n-2}W_m^{\frac{4}{n-2}}(\epsilon W_m+\phi)|\\
\\
\le C\left\{\begin{array}{l}
|\phi+\epsilon W_m|^{\frac{n+2}{n-2}},\quad\mbox{if}\quad |\phi|\ge W_m;\\
\\
W_m^{\frac{6-n}{n-2}}|\epsilon W_m+\phi|^2,\quad\mbox{if}\quad |\phi|\le W_m.
\end{array}\right.
\end{array}
$$

\medskip

By Lemma \ref{L8}, when $\lambda$ is large, if $|\phi|\ge W_m$, then $y\in \left(\cup_j (B_j\cap\Omega_j)\right)^c:=\Omega$. So
$$
\begin{array}{l}
\int_{{\mathbb R}^n}K_{\lambda}(y)(\bar{W}_m+\phi)^{\frac{n+2}{n-2}}Z_{i,j}dy\\
\\
=\int_{{\mathbb R}^n}K_{\lambda}(y)\left(W_m^{\frac{n+2}{n-2}}+\frac{n+2}{n-2}W_m^{\frac{4}{n-2}}(\epsilon W_m+\phi)\right)Z_{i,j}dy\\
\\
+O(1)\int_{\Omega}|\phi+\epsilon W_m|^{\frac{n+2}{n-2}}|Z_{i,j}|+O(1)\int_{{\mathbb R}^n}W_m^{\frac{6-n}{n-2}}|\epsilon W_m+\phi|^2|Z_{i,j}|.
\end{array}
$$

By Lemma \ref{L1}, Lemma \ref{L2} and similar argument as in the proof of Lemma \ref{L3}, we get easily that

$$
\begin{array}{l}
\int_{\Omega}|\phi+\epsilon W_m|^{\frac{n+2}{n-2}}|Z_{i,j}|\le C\int_{\Omega}|\phi|^{\frac{n+2}{n-2}}|Z_{i,j}|dy\\
\\
\le \frac{C\|\phi\|_*^{\frac{n+2}{n-2}}}{(\lambda l)^{\frac{4}{n-2}k}}\int_{\Omega}\sum_j\frac{1}{(1+|z-X^j|)^{\frac{n+2}{2}+\frac{n+2}{n-2}\tau-\frac{4}{n-2}k}}\frac{1}{(1+|z-X^i|)^{n-1}}dz\\
\\
\le \frac{C\|\phi\|_{*}^{\frac{n+2}{n-2}}}{(\lambda l)^{\frac{n}{2}+\tau\frac{n+2}{n-2}}},\quad \mbox{by Lemma \ref{L1}}.
\end{array}
$$

Using Lemma \ref{L1}, Lemma \ref{L2} and Lemma \ref{L5}, we can infer that
$$\begin{array}{l}
\int_{{\mathbb R}^n}W_m^{\frac{6-n}{n-2}}|\epsilon W_m+\phi|^2 |Z_{i,j}|dy\le \\
\\
C\int_{{\mathbb R}^n}\left(|\epsilon|^2W_m^{\frac{n+2}{n-2}}+W_m^{\frac{6-n}{n-2}}|\phi|^2\right)|Z_{i,j}|dy\\
\\
\le C\left(|\epsilon|^2+\left(\frac{\|\phi\|_*}{\lambda^{\tau-1}}\right)^2\right).
\end{array}
$$

\medskip

$$
|W_m^{\frac{4}{n-2}}-\sigma_i^{\frac{4}{n-2}}|\le C\left\{\begin{array}{c}
\hat{W}_{m,i}^{\frac{4}{n-2}},\quad\mbox{if}\quad \hat{W}_{m,i}>\sigma_i;\\
\\
\sigma_i^{\frac{4}{n-2}-1}\hat{W}_{m,i},\quad\mbox{if}\quad \hat{W}_{m,i}\le \sigma_i.
\end{array}\right.
$$

By Lemma \ref{L5}, we know  that $\hat{W}_{m,i} \leq \sigma_i$  in $B_i$ (we may need to shrink the ball a little bit) and $\hat{W}_{m,i}\geq \sigma_i$  in each $B_j$ with $j \neq i$.  Since $\l\phi,Z_{i,j}\r=0$, we can get,
$$
\begin{array}{l}
\int_{{\mathbb R}^n}K_{\lambda}(y)W_m^{\frac{4}{n-2}}\phi Z_{i,j}dy=\int_{{\mathbb R}^n}K_{\lambda}(y)\sigma_i^{\frac{4}{n-2}}\phi Z_{i,j}dy+o(1)\frac{1}{\lambda^{n}}\\
\\
=\int_{{\mathbb R}^n}(K_{\lambda}(y)-1)\sigma_i^{\frac{4}{n-2}}\phi Z_{i,j}dy+o(\frac{1}{\lambda^{n}})\\
\\
\le \frac{C}{\lambda^{\beta+1}}+o(\frac{1}{\lambda^{n}}).
\end{array}
$$

Similarly
$$
\begin{array}{l}
\int_{{\mathbb R}^n}K_{\lambda}(y)W_m^{\frac{4}{n-2}}\epsilon W_m Z_{i,j}dy=\int_{{\mathbb R}^n}(K_{\lambda}(y)-1)\epsilon_i\sigma_i^{\frac{n+2}{n-2}}Z_{i,j}dy\\
\\
+\frac{C|\epsilon|}{(\lambda l)^{n-2}}\le o(\frac{1}{\lambda^n}).
\end{array}
$$

For $n \geq 5$ it holds that
$$
\begin{array}{l}
|W_m^{\frac{n+2}{n-2}}-\sigma_i^{\frac{n+2}{n-2}}-\frac{n+2}{n-2}\sigma_i^{\frac{4}{n-2}}\hat{W}_{m,i}|\\
\\
\le C\left\{
\begin{array}{l}
\hat{W}_{m,i}^{\frac{n+2}{n-2}},\quad\mbox{if} \quad \hat{W}_{m,i}\ge \sigma_i;\\
\\
\hat{W}_{m,i}^2\sigma_i^{\frac{4}{n-2}-1}\le \hat{W}_{m,i}^{\frac{n}{n-2}}\sigma_i^{\frac{n}{n-2}-1} ,\quad\mbox{when}\quad  \hat{W}_{m,i}\le \sigma_i.
\end{array}\right.
\end{array}
$$

Using Lemma \ref{L1}, we deduce that
$$
\int_{{\mathbb R}^n}|K_{\lambda}(y)\hat{W}_{m,i}^{\frac{n+2}{n-2}}Z_{i,j}|dy\le \frac{C}{(\lambda l)^{n-1}},
$$
and
$$
\int_{{\mathbb R}^n}|K_{\lambda}(y)\hat{W}_{m,i}^{\frac{n}{n-2}}\sigma_i^{\frac{n}{n-1}-1}Z_{i,j}|dy\le \frac{C}{(\lambda l)^{n-1}}.
$$

For $s\neq i$, by change of variable for $s\neq i$,
$$
\begin{array}{l}
\frac{n+2}{n-2}\int_{{\mathbb R}^n}K_{\lambda}(y)\sigma_i^{\frac{4}{n-2}}\sigma_sZ_{i,j}dy=\frac{\partial}{\partial {P^i}_j}\int_{{\mathbb R}^n}K_{\lambda}(y)\sigma_i^{\frac{n+2}{n-2}}\sigma_sdy\\
\\
=\frac{\partial}{\partial {P^i}_j}\int_{{\mathbb R}^n}K(\frac{t+P^i}{\lambda})\sigma_{0,\Lambda_i}^{\frac{n+2}{n-2}}\sigma_{P^s-P^i,\Lambda_s}dt\\
\\
=\frac{1}{\lambda}\int_{{\mathbb R}^n}\frac{\partial K(\frac{t+P^i}{\lambda})}{\partial t_j}\sigma_{0,\Lambda_i}^{\frac{n+2}{n-2}}\sigma_{P^s-P^i,\Lambda_s}dt-\int_{{\mathbb R}^n}K(\frac{t+P^i}{\lambda})\sigma_{0,\Lambda_i}^{\frac{n+2}{n-2}}\frac{\partial\sigma_{P^s-P^i,\Lambda_s}}{\partial {P^s}_j}dt.
\end{array}
$$

From the above, using Lemma \ref{L1}, it is easy to see that
$$
|\int_{{\mathbb R}^n}K_{\lambda}(y)\sigma_i^{\frac{4}{n-2}}\hat{W}_{m,i} Z_{i,j}dy|\le \frac{C}{\lambda (\lambda l)^{n-2}}+\frac{C}{(\lambda l)^{n-1}}.
$$

Similarly
$$
|\sum_{s\neq i}(1+\epsilon_s)\l\frac{\partial\sigma_i}{\partial {P^i}_j},\sigma_s\r|\le \frac{C(1+|\epsilon|)}{(\lambda l)^{n-1}}.
$$

The above estimates  give the conclusion.
\end{proof}
\renewcommand{\theLemma}{A.12}
\begin{Lemma}\label{La.4} For some constant $C>0$ independent of $i$, $j$ and $m$,
$$
|\int_{{\mathbb R}^n}K_{\lambda}(x)\sigma_i^{\frac{n+2}{n-2}}\frac{\p\sigma_j}{\p \Lambda_j}dx-\int_{{\mathbb R}^n}\sigma_i^{\frac{n+2}{n-2}}\frac{\p\sigma_j}{\p \Lambda_j}|\le\frac{C}{|P^i-P^j|^{n-2}\lambda^2}.
$$
\end{Lemma}

\medskip

\begin{proof}
Take a $\delta>0$ small, such that $|K(x)-1|\le c|x|^{\beta}$ for some $c>0$ and $x\in B_{\delta}(0)$. Since $K(X^i)=1$, we get
\begin{equation*}
\begin{array}{c}
\int_{{\mathbb R}^n}K_{\lambda}(x)\sigma_i^{\frac{n+2}{n-2}}\frac{\p\sigma_j}{\p \Lambda_j}dx=\int _{{\mathbb R}^n}\sigma_i^{\frac{n+2}{n-2}}\frac{\p\sigma_j}{\p \Lambda_j}+\\
\\
\int_{B_{\delta\lambda}(X^i)\cup B_{\delta\lambda}(X^j)\cup(B^c_{\delta\lambda}(X^i)\cap B^c_{\delta\lambda}(X^j))}(K_{\lambda}(x)-1)\sigma_i^{\frac{n+2}{n-2}}\frac{\p\sigma_j}{\p \Lambda_j}dx,
\end{array}
\end{equation*}
and also
$$
\begin{array}{l}
|\int_{B_{\delta\lambda}(X^i)}(K_{\lambda}(x)-1)\sigma_i^{\frac{n+2}{n-2}}\frac{\p\sigma_j}{\p \Lambda_j}dx|\le C\int_{B_{\delta\lambda}(X^i)}\frac{|x-X^i|^{\beta}}{\lambda^{\beta}}\sigma_i^{\frac{n+2}{n-2}}\sigma_jdx\\
\\
\le C\int_{B_{\delta\lambda}(0)}(\frac{|x|^{\beta}}{\lambda^{\beta}}+\frac{|P^i-X^i|^{\beta}}{\lambda^{\beta}})\frac{1}{(1+|x|^2)^{\frac{n+2}{2}}}\frac{1}{(1+|x-P^j+P^i|^2)^{\frac{n-2}{2}}}dx\\
\\
\le \frac{C}{|P^i-P^j|^{n-2}}\int_{B_{\delta\lambda}(0)}(\frac{|x|^{\beta}}{\lambda^{\beta}}+\frac{|P^i-X^i|^{\beta}}{\lambda^{\beta}})\frac{1}{(1+|x|^2)^{\frac{n+2}{2}}}dx\\
\\
\le \frac{C}{|P^i-P^j|^{n-2}}(\frac{1}{\lambda^2}+\frac{1}{\lambda^{\beta}})\le \frac{C}{|P^i-P^j|^{n-2}\lambda^2}.
\end{array}
$$

Similarly,
$$
\begin{array}{l}
|\int_{B_{\delta\lambda}(X^j)}(K_{\lambda}(x)-1)\sigma_i^{\frac{n+2}{n-2}}\frac{\p\sigma_j}{\p \Lambda_j}dx|\\
\\
\le \frac{C}{|P^i-P^j|^{n+2}}(\lambda^2+\lambda^{3-\beta})\le \frac{C}{|P^i-P^j|^{n-2}\lambda^2}.
\end{array}
$$

Lastly,
\begin{equation*}\label{La.4.1}\tag{A13}
\begin{array}{l}
|\int_{B^c_{\delta\lambda}(X^i)\cap B^c_{\delta\lambda}(X^j)}(K_{\lambda}(x)-1)\sigma_i^{\frac{n+2}{n-2}}\frac{\p\sigma_j}{\p \Lambda_j}dx|\\
\\
\le C\int_{B^c_{\delta\lambda}(X^i)\cap B^c_{\delta\lambda}(X^j)}\sigma_i^{\frac{n+2}{n-2}}\sigma_jdx\\
\\
\le C\int_{B^c_{\delta\lambda}(0)\cap B^c_{\delta\lambda}(X^j-X^i)}\frac{1}{(1+|x|^2)^{\frac{n+2}{2}}(1+|x-P^j+P^i|^2)^{\frac{n-2}{2}}}dx\\
\end{array}
\end{equation*}

Let $z=P^j-P^i$ and $2d=|z|$. To estimate (\ref{La.4.1}), we will use the method used by Wei-Yan in \cite{WY}. Since $d>\delta\lambda$ when $l$ is large, we can split
$$
B^c_{\delta\lambda}(0)\cap B^c_{\delta\lambda}(X^j-X^i)=A_1\cup A_2\cup A_3,$$
where $A_1=B_d(0)\setminus B_{\delta\lambda}(0)$, $A_2=B_d(X^j-X^i)\setminus B_{\delta\lambda}(X^j-X^i)$ and $A_3=B^c_d(0)\cap B^c_d(X^j-X^i)$.

$$
\begin{array}{l}
\int_{A_1}\frac{1}{(1+|x|^2)^{\frac{n+2}{2}}(1+|x-z|^2)^{\frac{n-2}{2}}}dx\le \frac{C}{d^{n-2}}\int_{A_1}\frac{1}{(1+|x|^2)^{\frac{n+2}{2}}}dx\\
\\
\le \frac{C}{|P^i-P^j|^{n-2}\lambda^2}.
\end{array}
$$

Similarly it holds that
$$
\int_{A_2}\frac{1}{(1+|x|^2)^{\frac{n+2}{2}}(1+|x-z|^2)^{\frac{n-2}{2}}}dx\le \frac{C}{|P^i-P^j|^n}.
$$

On $A_3$, from \cite{WY},
$$
\frac{1}{(1+|x|^2)^{\frac{n+2}{2}}(1+|x-z|^2)^{\frac{n-2}{2}}}\le \frac{C}{|x|^{n-2}(1+|x|)^{n+2}}\le \frac{C}{|x|^{2n}},
$$therefore, we infer that
$$
\int_{A_3}\frac{1}{(1+|x|^2)^{\frac{n+2}{2}}(1+|x-z|^2)^{\frac{n-2}{2}}}dx\le\frac{C}{|P^i-P^j|^n},
$$
which gives

$$
|\int_{B^c_{\delta\lambda}(X^i)\cap B^c_{\delta\lambda}(X^j)}(K_{\lambda}(x)-1)\sigma_i^{\frac{n+2}{n-2}}\frac{\p\sigma_j}{\p \Lambda_j}dx|\le \frac{C}{|P^i-P^j|^{n-2}\lambda^2}.
$$

The conclusion of Lemma follows easily.
\end{proof}

\renewcommand{\theLemma}{A.13}
\begin{Lemma}\label{La.5} For some constant $C>0$ independent of $i$, $j$, $m$,
$$
|\int_{{\mathbb R}^n}K_{\lambda}(x)\sigma_j^{\frac{4}{n-2}}\sigma_i\frac{\p\sigma_j}{\p\Lambda_j}dx-\frac{n-2}{n+2}\int_{{\mathbb R}^n}\sigma_i^{\frac{n+2}{n-2}}\frac{\p\sigma_j}{\p\Lambda_j}dx|\le \frac{C}{|P^i-P^j|^{n-2}\lambda^2}.
$$
\end{Lemma}
\begin{proof}
$$
\begin{array}{l}
\int_{{\mathbb R}^n}K_{\lambda}(x)\sigma_j^{\frac{4}{n-2}}\sigma_i\frac{\p\sigma_j}{\p\Lambda_j}dx=\frac{n-2}{n+2}\int_{{\mathbb R}^n}K_{\lambda}(x)\frac{\p\sigma_j^{\frac{n+2}{n-2}}}{\p\Lambda_j}\sigma_idx\\
\\
=\frac{n-2}{n+2}\int_{{\mathbb R}^n}\frac{\p\sigma_j^{\frac{n+2}{n-2}}}{\p\Lambda_j}\sigma_idx+\frac{n-2}{n+2}\int_{{\mathbb R}^n}(K_{\lambda}(x)-1)\frac{\p\sigma_j^{\frac{n+2}{n-2}}}{\p\Lambda_j}\sigma_idx\\
\\
=\frac{n-2}{n+2}\int_{{\mathbb R}^n}\sigma_i^{\frac{n+2}{n-2}}\frac{\p\sigma_j}{\p\Lambda_j}dx+\frac{n-2}{n+2}\int_{{\mathbb R}^n}(K_{\lambda}(x)-1)\frac{\p\sigma_j^{\frac{n+2}{n-2}}}{\p\Lambda_j}\sigma_idx.
\end{array}
$$

The second error term in the above can be similarly estimated  as in Lemma \ref{La.4} and  the proof is thus completed.
\end{proof}

\bigskip

\section{Appendix B}

We recall some results proved in \cite{L2}, in a form convenient for our application.

Let $\{K_i\}$ be a sequence of functions satisfying, for some constant $A\ge 1$,
\begin{equation*}\tag{B1}
\frac{1}{A}\le K_i\le A,\quad\mbox{in}\quad B_2\quad \forall i.
\end{equation*}
where $B_2$ is the ball in ${\mathbb R}^n$ of radius 2 centered at the origin.

For $\beta>n-2$, recall that $\{K_i\}$ satisfies $(*)_{\beta}$ for some positive constants $L_1$ and $L_2$ (independent of $i$) in $B_2$ if $\{K_i\}\subset C^{[\beta]-1,1}(B_2)$ satisfies
    $$|\nabla K_i|\le L_1,\quad \mbox{in}\quad B_2,$$ and, if $\beta\ge 2$, that
    $$
    |\nabla^s K_i(y)|\le L_2|\nabla K_i(y)|^{\frac{\beta-s}{\beta-1}},\quad\mbox{for all}\quad  2\le s\le [\beta],\quad y\in B_2.$$

\begin{remark}
\em
Conditions {\bf (H1)}, {\bf (H2)} and {\bf (H3)} guarantee that $K$ satisfies $(*)_{\beta}$ in a neighborhood of 0.
\end{remark}

\renewcommand{\theProposition}{B.1}
\begin{Proposition}\label{p3}
For $n\ge 3$, $A\ge 1$, $L_1$, $L_2>0$ and $\beta>n-2$, let $\{K_i\}$ be a sequence of functions satisfying (B1) and $(*)_{\beta}$ for $L_1$ and $L_2$. Let $\{u_i\}$ be a sequence of $C^2$ solutions of

$$
-\Delta u_i=K_iu_i^{\frac{n+2}{n-2}},\quad u_i>0,\quad B_2,
$$satisfying, for some constant $a$ independent of $i$,
$$
\|u_i\|_{L^{\frac{2n}{n-2}}(B_2)}\le a<\infty,\quad\forall i.
$$ Then, after passing to a subsequence, $\{u_i\}$ either stays bounded in $B_1$, or has only isolated simple blow up points in $B_1$.
\end{Proposition}

See Definition 0.3 in \cite{L2} for the definition of isolated simple blow up points.
\renewcommand{\theProposition}{\arabic{Proposition}}
\begin{proof}
It is easy to see, using the equation of $u_i$, that
$$
\|\nabla u_i\|_{L^2(B_{\frac{3}{2}})}\le C(n,A,a):=C_a.
$$
Let $\delta_0>0$ be the small constant given in Proposition 2.1 in \cite{L}, and fix a positive integer $k$ such that
$$
C^2_a+a^{\frac{2n}{n-2}}\le \delta_0k,
$$
For $r_l=1+\frac{l-1}{2k}$, $1\le l\le k+1$, let
$$
A_l=\{x|r_l\le|x|\le r_{l+1}\},\quad 1\le l\le k.
$$

Since
$$\begin{array}{l}
\sum_{l=1}^k\int_{A_l}(|\nabla u_i|^2+u_i^{\frac{n+2}{n-2}})\le\int_{B_{\frac{3}{2}}}(|\nabla u_i|^2+u_i^{\frac{n+2}{n-2}})\le C_a^2+a^{\frac{2n}{n-2}}\le\delta_0k,
\end{array}
$$there exist some $1\le l\le k$ and a subsequence of $\{u_i\}$ (still denoted as $\{u_i\}$) such that

$$
\int_{A_l}(|\nabla u_i|^2+u_i^{\frac{n+2}{n-2}})\le\delta_0, \quad\forall i.
$$

It follows from Proposition 2.1 in \cite{L} that
$$
\|u_i\|_{L^{\infty}(\hat{A}_l)}\le C(\delta_0,n, A, a),
$$
where $\hat{A}_l=\{x|r_l+\frac{1}{8k}<|x|<r_{l+1}-\frac{1}{8k}\}$. Using this estimate,  we work on the ball $B_{r_{l+1}-\frac{1}{8k}}$. Then the proofs of Proposition 4.1, Proposition 4.2 and Theorem 4.2 in \cite{L2} apply, in view of the fact that $\{u_i\}$ stays bounded in the shell $\hat{A}_l$. We obtain the conclusion.
\end{proof}

\end{document}